\documentclass[reqno,11pt]{amsart}
\usepackage[english]{babel}
\usepackage{amsmath,amsfonts,amsthm,amssymb,amsbsy,upref,color,graphicx,amscd,hyperref,enumerate,bbm,comment, algorithm, algorithmic}
\usepackage{mathtools,url}
\usepackage[a4paper,margin=2.8truecm]{geometry}
\usepackage[active]{srcltx}

\DeclareMathOperator{\dive}{div}

\DeclareMathOperator{\spt}{spt}

\DeclareMathOperator{\dist}{dist}

\def\ds{\displaystyle}
\def\eps{{\varepsilon}}
\def\N{\mathbb{N}}
\def\O{\Omega}

\def\R{\mathbb{R}}

\def\EE{\mathcal{E}}

\def\HH{\mathcal{H}}

\def\M{\mathcal{M}}

\def\V{\mathcal{V}}
\def\symd{{\scriptstyle\Delta}}
\def\ecart{\noalign{\medskip}}
\def\ep{{\varepsilon}}

\newcommand{\be}{\begin{equation}}
\newcommand{\ee}{\end{equation}}
\newcommand{\bib}[4]{\bibitem{#1}{\sc#2: }{\it#3. }{#4.}}
\newcommand{\cp}{\mathop{\rm cap}\nolimits}

\numberwithin{equation}{section}
\theoremstyle{plain}
\newtheorem{theo}{Theorem}[section]
\newtheorem{lemm}[theo]{Lemma}

\newtheorem{prop}[theo]{Proposition}
\newtheorem{defi}[theo]{Definition}

\theoremstyle{remark}
\newtheorem{rema}[theo]{Remark}

\title[Optimal potentials on unbounded domains]{On the existence of optimal potentials on unbounded domains}

\author[G. Buttazzo]{Giuseppe Buttazzo}
\address{Dipartimento di Matematica, Universit\`a di Pisa, Largo B. Pontecorvo 5, 56126 Pisa, ITALY}
\email{giuseppe.buttazzo@unipi.it}

\author[J. Casado Diaz]{Juan Casado-D\'{\i}az}
\address{Dpto. Ecuaciones Diferenciales y An\'alisis Num\'erico, Universidad de Sevilla, C/ Tarfia s/n. Aptdo 1160, 41080 Sevilla, SPAIN}
\email{jcasadod@us.es}

\author[F. Maestre]{Faustino Maestre}
\address{Dpto. Ecuaciones Diferenciales y An\'alisis Num\'erico, Universidad de Sevilla, C/ Tarfia s/n. Aptdo 1160, 41080 Sevilla, SPAIN}
\email{fmaestre@us.es}

\date{\today}

\begin{document}

\begin{abstract}
We consider elliptic equations of Schr\"odinger type with a right-hand side fixed and with the linear part of order zero given by a potential $V$. The main goal is to study the optimization problem for an integral cost depending on the solution $u_V$, when $V$ varies in a suitable class of admissible potentials. These problems can be seen as the natural extension of shape optimization problems to the framework of potentials. The main result is an existence theorem for optimal potentials, and the main difficulty is to work in the whole Euclidean space $\R^d$, which implies a lack of compactness in several crucial points. In the last section we present some numerical simulations.
\end{abstract}

\maketitle

\textbf{Keywords:} Optimal potentials, Schr\"odinger operators, shape optimization, unbounded domains, capacitary measures.

\textbf{2010 Mathematics Subject Classification:} 49J45, 49Q10, 35J10, 49A22, 35J25, 49B25.

\section{Introduction}\label{sintr}

In the present paper we consider the optimization problem
\be\label{1.1}
\min\left\{\int j(x,u_V,\nabla u_V)\,dx\ :\ V\in\V\right\}.
\ee
The problem above can be seen as an optimal control problem where:
\begin{itemize}
\item the control variable $V$ is a nonnegative potential;
\item $u_V$ denotes the unique solution of the state equation, which is a PDE of Schr\"odinger type, formally written as
\be\label{E1.1}-\Delta u+V(x)u=f\hbox{ in }\R^d,\qquad u\in H^1(\R^d)\ee
with the right-hand side $f$ fixed;
\item the cost function $j(x,s,z)$ satisfies suitable mild conditions;
\item the class $\V$ of admissible controls is of the form
$$\V=\left\{\Psi(V)\le1\right\}$$
with $\Psi(V)$ an integral functional satisfying suitable conditions.
\end{itemize}
The precise assumptions on $j,f,\psi$ will be given in Section \ref{sprel}. Here we want to stress that the ambient space is the whole $\R^d$; working on the whole space $\R^d$ represents indeed the main difficulty, because several compactness theorems fail and parts of minimizing sequences $(u_n,V_n)$ may {\it``escape to infinity''}. We recall that similar problems on a bounded ambient space have been considered in \cite{bubuve}, \cite{bgrv14}, \cite{buve17}.

For simplicity, along all the paper, the notation of function spaces $L^2$, $H^1$ and similar, without the indication of the domain of definition, is used when the domain is the whole $\R^d$. Similarly, the absence of the domain of integration in an integral means that the integral is made on the whole $\R^d$.

Optimization problems of the form \eqref{1.1} are the natural extension to the class of potentials of shape optimization problems, that are written as
$$
\min\left\{\int j(x,u_\O,\nabla u_\O)\,dx\ :\ |\O|\le1\right\},
$$
where $u_\O$ denotes the unique solution of the Dirichlet problem
\be\label{1.2}-\Delta u=f\hbox{ in }\O,\qquad u\in H^1_0(\O).\ee
In fact, a domain $\O$ can be represented by the potential formally written as
$$V(x)=\begin{cases}
0&\hbox{if }x\in\O\\
+\infty&\hbox{if }x\notin\O.
\end{cases}$$
in the sense that the PDE \eqref{E1.1} becomes the PDE \eqref{1.2}. To be rigorous, when $\mu$ is a capacitary measure (see Section \ref{sprel} for the rigorous definition) the PDE formally written as
$$-\Delta u +\mu u=f,\qquad u\in H^1$$
has to be intended in the weak form as $u\in H^1\cap L^2_\mu$ and
$$\int\nabla u\cdot\nabla v\,dx+\int uv\,d\mu=\int fv\,dx\qquad\forall v\in H^1\cap L^2_\mu.$$
Similarly, a capacitary measure $\mu$ can be decomposed as
$$\mu=\mu^a+\mu^s+\mu^\infty$$
where $\mu^a$ and $\mu^s$ are respectively the absolutely continuous and the singular parts of $\mu$ with respect to the Lebesgue measure, and $\mu^\infty$ is the infinite part, of the form
$$\mu^\infty(E)=\begin{cases}
0&\hbox{if }\cp(E\cap K)=0\\
+\infty&\hbox{if }\cp(E\cap K)>0
\end{cases}$$
for some quasi-closed set $K$. Then, the class $\V$ of admissible potentials has to be intended in the sense of integral functionals over measures as
$$\Psi(\mu)=\int\psi(\mu^a)\,dx+C_\psi\mu^s(\R^d)+\psi(\infty)\cp(K)$$
where $\psi$ is a given nonnegative convex function and
$$\psi(\infty)=\lim_{t\to+\infty}\psi(t),\qquad C_\psi=\lim_{t\to+\infty}\frac{\psi(t)}{t}\;.$$
The main result of the paper is an existence theorem (Theorem \ref{ThExsoPo}) for minimizers of problem \eqref{1.1}. The detailed presentation of the optimization problem is given in Section \ref{sprel}; Section \ref{sproo} contains the proofs of the results, while in Sections \ref{snece} and \ref{sprooCO} we collected some necessary conditions of optimality for the solutions of problem \eqref{1.1}. Finally, in Section \ref{snume} we present some numerical simulations.
\section{Preliminaries and statement of the main results}\label{sprel}

\subsection{Preliminaries about capacity}\label{sscapa}

In the paper we use the key notion of capacity; for the sake of completeness we recall here its definition together with the terminology we adopt; the reader interested in this topic can find details and proofs on the facts below on \cite{bubu05}.

For a subset $E\subset\R^d$ its {\it capacity} is defined by
$$\cp(E)=\inf\left\{\int|\nabla u|^2\,dx+\int u^2\,dx\ :\ u\in H^1,\ u\ge 1\ \hbox{in a neighborhood of } E\right\}.$$
If a property $P(x)$ holds for all $x\in\R^d$, except for the elements of a set $E$ of capacity zero, we 
say that $P(x)$ holds {\it quasi-everywhere} (shortly {\it q.e.}), whereas the expression {\it almost
everywhere} (shortly {\it a.e.}) refers, as usual, to the Lebesgue measure, which we often denote by $|\cdot|$.

\begin{defi}
A set $\O\subset\R^d$ is called {\it quasi-open} (respectively {\it quasi-closed}) if there exists a sequence $(A_n)$ of open (respectively closed) sets such that $\cp(\O\symd A_n)$ tends to zero, where $\symd$ denotes the symmetric difference of sets.
\end{defi}

\begin{rema}
It can be seen that a set $\O$ is quasi-open set if and only if $\O=\{u>0\}$ for a suitable function $u\in H^1_{loc}$ and similarly that $K$ is quasi-closed if and only if $K=\{u\ge0\}$ for a suitable function $u\in H^1_{loc}$.
\end{rema}

\begin{defi}\label{capm}
A nonnegative Borel measure $\mu$ on $\R^d$ (possibly taking the value $+\infty$) is called a {\it capacitary measure} if for every Borel set $E\subset\R^d$ we have
\[\begin{cases}
\cp(E)=0\ \Longrightarrow\ \mu(E)=0;\\
\mu(E)=\inf\big\{\mu(\O)\ :\ \O\supset E,\ \O\hbox{ quasi-open}\big\}.
\end{cases}\]
We denote by $\M_{cap}$ the class of all capacitary measures on $\R^d$, while $\M$ stands for the class of all nonnegative Radon measures on $\R^d$.
\end{defi}

\begin{rema}
Since Sobolev functions $u\in H^1_{loc}$ are defined up to a capacity null set, the quantity $\int u^2\,d\mu$ is well defined for every $u\in H^1_{loc}$ and for every capacitary measure $\mu$.
\end{rema}

\begin{rema} \label{remEjMe}
Every measure $\mu$ which is absolutely continuous with respect to the Lebesgue measure, is a capacitary measure. Indeed, if $\cp(E)=0$, then $|E|=0$ and so $\mu(E)=0$. In particular, if $a(x)$ is a nonnegative function in $L^1_{loc}$, the measure $a(x)\,dx$ is of capacitary type. Analogously, if $S$ is a $d-1$ regular manifold and $a(x)$ is a nonnegative function in $L^1_{loc}(S)$ the measure $a(x)\,d\HH^{d-1}\lfloor S$ is of a capacitary type, where $\HH^{d-1}$ denotes the $d-1$ Hausdorff measure. Finally, for every quasi-closed set $K\subset\R^d$ the measure
\be\label{defmiE}\infty_K(E)=\begin{cases}
0&\hbox{if }\cp(E\cap K)=0\\
+\infty&\hbox{if }\cp(E\cap K)>0
\end{cases}\ee
is a capacitary measure.
\end{rema}

\subsection{The $\gamma$-convergence}\label{ssgamm}

We recall the definition of $\Gamma$-convergence.

\begin{defi}\label{deGco}
Given a metric space $X$, we say that a sequence of functionals $F_n:X\to(-\infty,\infty]$ $\Gamma$-converges to a functional $F:X\to (-\infty,\infty]$ if
\[\begin{split}
&\forall\,u\in X,\ \forall\,u_n\to u\hbox{ in }X\quad F(u)\le\liminf_{n\to\infty}F_n(u_n);\\
&\forall\,u\in X,\ \exists\,u_n\to u\hbox{ such that }\quad F(u)=\lim_{n\to\infty}F_n(u_n).
\end{split}\]
A sequence $u_n$ satisfying the equality above is said to be a recovery sequence for $u$.
\end{defi}

\begin{rema} We recall that the main interest of the $\Gamma$-convergence is the study of the asymptotic behavior of the minimum points and values of the functionals $F_n$. Namely, it is well known that if $(u_n)$ is a compact sequence in $X$ such that $F_n(u_n)-\inf_XF_n$ tends to zero, then $F$ admits a minimum and every cluster point $u$ of $(u_n)$ in $X$ satisfies $F(u)=\min_XF$.
\end{rema} 
\begin{rema} \label{RemRS} In the case where $X$ is a vector space and the functionals $F_n$ in Definition \ref{deGco} are quadratic, i.e. $F_n(u)=a_n(u,u)$ with $a_n$ bilinear, it is known that $u_n$ is a recovery sequence for $u$ if and only if $u_n$ converges to $u$ in $X$ and satisfies
$$\lim_{n\to\infty} a_n(u_n,v_n)=0$$
for every $v_n\to0$ in $X$ with $\limsup_n F_n(v_n)<\infty$.
\end{rema}
In our case we are interested in the sequence of functionals 
$$u\to \int|\nabla u|^2\,dx+\int u^2\,d\mu_n$$
where $\mu_n$ is a sequence in $\M_{cap}$. In the case where we deal with functions defined in a bounded open set of $\R^d$ it has been proved in \cite{dmmo87} that the class of these functionals is closed for the $\Gamma$-convergence. In our case we deal with functions which are defined in the whole of $\R^d$. The first difficulty is to choose the good spaces of functions where these functionals are well defined, and the usual Sobolev space $H^1$ is not a good choice since Poincar\'e's inequality does not hold for functions in $H^1$.

\begin{defi} Denoting by $W:\R^d\to\R$ the function
$$W(x)={1\over1+|x|}\ \hbox{ if }d\neq2,\qquad W(x)={1\over(1+|x|)\log(2+|x|)}\ \hbox{ if }d=2,$$
we define $L$ as the space
\be\label{defiL}L=\left\{u:\R^d\to\R:\quad Wu\in L^2\right\},\ee
endowed with the norm
$$\|u\|_L=\|uW\|_{L^2},$$
and $H$ as
$$H=\big\{u\in L\cap H^1_{loc}:\ \nabla u\in (L^2)^d\big\}.$$
\end{defi}

\begin{prop}\label{DenCiH}
The space $C^\infty_c$ is dense in $H$. The usual norm in $H$, given by
$$\|u\|^2=\|u\|^2_L+\|\nabla u\|^2_{L^2}\;,$$
is equivalent to
\be\label{norHd>2}
\|u\|_H=\|\nabla u\|_{L^2}\qquad\hbox{if }d\ge3,
\ee
\be\label{norHd<=2}
\|u\|^2_H=\|u\|^2_{L^2(B(0,1))}+\|\nabla u\|^2_{L^2}\quad\hbox{if }d=1,2.
\ee
\end{prop}

\begin{rema}\label{carHd>2}
From the Sobolev embedding theorem, we also have, if $d\ge3$,
$$H=\big\{u\in L^{2d/(d-2)}\ :\ \nabla u\in (L^2)^d\big\}.$$
\end{rema}

\begin{defi} We say that a sequence $(\mu_n)$ in $\M_{cap}$ $\gamma$-converges to a measure $\mu\in\M_{cap}$ if the sequence of functionals $F_n:H\to [0,\infty]$ defined by
$$F_n(u)=\int|\nabla u|^2\,dx+\int u^2\,d\mu_n$$
$\Gamma$-converges in $H$, endowed with the topology of $L^2_{loc}$, to the functional $F$ given by
$$F(u)=\int|\nabla u|^2\,dx+\int u^2\,d\mu.$$
\end{defi}

When $(\mu_n)$ are defined in a bounded open set $\O$ of $\R^d$ (in the sense that $\mu_n=+\infty$ outside $\O$), it has been proved in \cite{dmmo87} (see also \cite{bdm91}, \cite{bdm93}, \cite{DMGA}) that every sequence of measures in $\M_{cap}$ contains a subsequence which $\gamma$-converges. The following theorem extends this result to measures defined in the whole of $\R^d$.

\begin{theo}\label{gacomp}
For every sequence $\mu_n\in\M_{cap}$ there exists a subsequence (still denoted by $\mu_n$) and $\mu\in\M_{cap}$ such that $\mu_n$ $\gamma$-converges to $\mu$.
\end{theo}

Using Theorem \ref{gacomp} we can prove the following proposition which is useful to study the asymptotic behavior of the solutions of the elliptic problems
$$-\Delta u_n+u_n\mu_n=f\ \hbox{ in }\R^d$$
when $d=1,2$.

\begin{prop}\label{DePd12}
Assume $d=1,2$ and let $\mathcal{U}\subset\M_{cap}$ be a subset which is closed for the $\gamma$-convergence and which does not contain the null measure. Then, there exists $C>0$ such that
$$\|u\|_L\le C\big(\|\nabla u\|_{L^2}+\|u\|_{L^2_\mu}\big),\qquad\forall\mu\in\mathcal{U},\ \forall u\in H\cap L^2_\mu.$$
\end{prop}

\begin{rema}\label{DefSoEDP} Proposition \ref{DePd12} in the case $d=1,2$ and \eqref{norHd>2} if $d\geq 3$ combined with the Lax-Milgram theorem allows to prove the existence and uniqueness of solutions for the problem
\be\label{pbelimu}
-\Delta u+\mu u=f\ \hbox{ in }\R^d,\qquad u\in H\cap L^2_\mu,
\ee
where $f$ belongs to $H'$ and $\mu$ to $\M_{cap}$, with $\mu$ not the null measure if $d=1,2$. We recall that since $\mu$ is not in general a Radon measure, equation \eqref{pbelimu} has not to be intended in the distributional sense. Namely, problem \eqref{pbelimu} has to be understood in the variational sense
$$\begin{cases}
u\in H\cap L^2_\mu\\
\ds\int\nabla u\cdot\nabla v\,dx+\int uv\,d\mu=\langle f,v\rangle\qquad \forall\,v\in H\cap L^2_\mu.
\end{cases}$$
Taking for instance $\mu$ as the measure $\infty_K$ given by \eqref{defmiE}, we get that equation \eqref{pbelimu} can be read as
$$-\Delta u=f\ \hbox{ in }\R^d\setminus K,\qquad u=0\ \hbox{ on }K.$$
Since the functions with compact support are dense in $H\cap L^2_\mu$, the condition $u\in H\cap L^2_\mu$ can be interpreted as ``$u=0$ at infinity''. Observe however that for $d=1,2$ there exist elements in $H$ which tend to infinity at infinity. We also recall that in the distributional sense every $f\in H'$ can be written as
$$f=f_1-\dive f_2$$
with $f_1W^{-1}\in L^2$ and $f_2\in (L^2)^d$.
\end{rema}

Using the $\gamma$-convergence we can now prove the following result about the asymptotic behavior of the solution of \eqref{pbelimu} when $\mu$ varies.

\begin{prop}\label{ConsoEc}
Let $\mu_n\in\M_{cap}$ be a sequence which $\gamma$-converges to $\mu$, where for $d=1,2$ the measures $\mu_n$ and the measure $\mu$ are not the null measure. Then, for every $f\in H'$, the solutions $u_n$ of
\be\label{pbun}
-\Delta u_n+\mu_n u_n=f\ \hbox{ in }\R^d\qquad u_n\in H\cap L^2_\mu,
\ee
satisfy
\be\label{convun}
u_n\rightharpoonup u\ \hbox{ in }H,\qquad u_n\rightarrow u\ \hbox{ in }W^{1,p}_{loc}\quad\hbox{for every }p<2,
\ee
where $u$ is the solution of \eqref{pbelimu}.
\end{prop}

\subsection{The optimization problem}\label{ssprob}

The optimization problems we aim to consider are written in the form
\be\label{optpb}
\min\left\{\int j(x,u,\nabla u)\,dx\ :\ -\Delta u+\mu u=f,\ \Psi(\mu)\le1,\ \mu\ge\nu\right\},
\ee
where $f\in H'$, $\nu\in\M_{cap}$, the function $j(x,s,\xi)$ verifies some suitable conditions, and the functional $\Psi$ is suitably defined. The meaning of the PDE $-\Delta u+\mu u=f$ which is posed in the whole space $\R^d$ has been explained in Remark \ref{DefSoEDP}.

In several cases (always if $d\ge3$) the measure $\nu$ can be chosen as the null measure and then condition $\mu\ge\nu$ is not a restriction. Concerning the function $j:\R^d\times\R\times\R^d\to\R\cup\{+\infty\}$ we assume (see for instance \cite{busc}) that $j(x,s,\xi)$ is measurable in $x$, lower semicontinuous in $(s,\xi)$ and verifies the inequality
\be\label{hypoj}
-g_1(x)|s|-g_2(x)|s|^2-h_1(x)|\xi|-h_2(x)|\xi|^q\le j(x,s,\xi)
\ee
for a.e. $x\in\R^d$ and for all $(s,\xi)\in\R\times\R^d$, where $1<q<2$, $g_1,g_2,h_1,h_2\ge0$ a.e. in $\R^d$, $g_1W^{-1}\in L^2$, $g_2\in L^\infty$, $h_1\in L^2$, $h_2\in L^{2/(2-q)}$ and
\be\label{hypoab}
\begin{cases}
\ds{\rm ess\,}\lim_{\hskip-10pt x\to\infty} g_2(x)|x|^2=0&\hbox{if }d\neq2,\\
\ds{\rm ess\,}\lim_{\hskip-10pt x\to\infty} g_2(x)|x|^2\log^2|x|=0&\hbox{if }d=2.
\end{cases}
\ee
The functional $\Psi$ acts on capacitary measures; every $\mu\in\M_{cap}$ can be uniquely written in the form
\be\label{decompmu}\mu=\mu^a+\mu^s+\mu^\infty
\ee
where $\mu^\infty$ is of the form $\infty_K$ for some quasi-closed set $K$, $\mu^a$ is absolutely continuous with respect to the Lebesgue measure, $\mu^s\in\M_{cap}$ is singular with respect to the Lebesgue measure, $\mu^a(K)=\mu^s(K)=0$, and $\mu^a+\mu^s$ is $\sigma$-finite. 

We consider $\Psi$ of the form
\be\label{defpsi}
\Psi(\mu)=\int\psi(\mu^a+\mu^s)+\psi(\infty)\cp(K)
\ee
where $\psi:\R^+\to[0,+\infty]$ is a convex and lower semicontinuous function and
$$\psi(\infty)=\lim_{t\to+\infty}\psi(t).$$
The integral above is intended in the sense of convex integral on measures, more precisely as
$$\int\psi(\mu^a+\mu^s)=\int\psi(\mu^a)\,dx+C_\psi\mu^s(\R^d),$$
where
\be\label{kpsi}
C_\psi=\lim_{t\to+\infty}\frac{\psi(t)}{t}\;.
\ee

\begin{rema}\label{Rea0}
It is known that $C_\psi$ defined by \eqref{kpsi} always exists and satisfies
$$C_\psi=\sup\big\{\tau\in\partial\psi(s)\ :\ s\in (0,\infty)\big\},$$
where $\partial\psi(s)$ denotes the subdifferential of the convex function $\psi$ at $s$. If $C_\psi=0$, we have that $\psi$ is a decreasing function and therefore there exists a finite limit at infinity of $\psi$. This proves that
$$\lim_{s\to+\infty}\frac{\psi(s)}{s}=0\ \iff\ \exists\lim_{s\to\infty}\psi(s)=\ell<+\infty.$$
If $\ell>0$ then $\psi(s)\ge\ell$ for every $s\in\R$ and so $\Psi$ is the trivial functional defined by
$$\Psi(\mu)=\infty\qquad\forall\,\mu\in\M_{cap}.$$
So, in the case $C_\psi=0$, we assume that
\be\label{supa0}
\psi(\infty)=\lim_{s\to+\infty}\psi(s)=0,
\ee
so that the functional $\Psi$ reduces to
$$\Psi(\mu)=\int\psi(\mu^a)\,dx\qquad\forall\mu\in\M_{cap}.$$
If $C_\psi>0$, the function $\psi$ attains a minimum. As above, in order to have $\Psi$ not trivial, we assume
$$\min_{s\in[0,\infty)}\psi(s)=0.$$
In addition, $C_\psi>0$ implies
$$\psi(\infty)=\lim_{s\to+\infty}\psi(s)=+\infty,$$
so that the functional $\Psi$ reduces to
$$\Psi(\mu)=\begin{cases}
\int\psi(\mu^a)\,dx+C_\psi\mu^s(\R^d)&\hbox{if }\mu^\infty=0,\\
+\infty&\hbox{otherwise,}
\end{cases}$$
which in the case $C_\psi=+\infty$ just gives
$$\Psi(\mu)=\begin{cases}
\int\psi(\mu^a)\,dx&\hbox{if }\mu^s=\mu^\infty=0,\\
+\infty&\hbox{otherwise.}
\end{cases}$$
\end{rema}

\begin{prop}\label{psisci}
The functional $\Psi$ is lower semicontinuous on $\M_{cap}$ with respect to the $\gamma$-convergence.
\end{prop}

As a consequence of the above result we can now prove the following theorem.

\begin{theo}\label{ThExsoPo}
Assume $j:\R^d\times\R\times\R^d\to\R\cup\{+\infty\}$ measurable in $x$, lower semicontinuous in $(s,\xi)$ and satisfying
\eqref{hypoj} for some $q\in(1,2)$, $g_1,g_2,h_1,h_2\ge0$ a.e. in $\R^d$, $g_1W^{-1}\in L^2$, $g_2\in L^\infty$ satisfying \eqref{hypoab}, $h_1\in L^2$, $h_2\in L^{2/(2-q)}$. We consider a function $\psi:\R^+\to[0,\infty]$ convex and lower semicontinuous and a measure $\nu\in\M_{cap}$ such that there exists $\hat\mu\in\M_{cap}$ satisfying
$$\hat\mu\ge\nu,\qquad\Psi(\hat\mu)\le1,$$
with $\Psi$ defined by \eqref{defpsi}. Moreover, if $d=1,2$ we assume that:
\be\label{Conpsinu}
\hbox{either $\psi(0)>0$ or $\nu$ is not the null measure.}\ee
Then, for every $f\in H'$, problem \eqref{optpb} has at least one solution.
\end{theo}

\section{Proofs of the results of section \ref{sprel}}\label{sproo}

\begin{proof}[Proof of Proposition \ref{DenCiH}] Let us prove the density of $C^\infty_c$ in $H$. Taking into account the density of $C^\infty_c(\O)$ in $H^1_0(\O)$, for every bounded open set $\O\subset\R^d$, it is enough to show that for every $u\in H$, there exists a sequence $u_n\in H^1$ with compact support, which converges to $u$ in $H$.

We first consider the case $d\neq 2$. Then, we take $u_n=u\varphi_n$ with $\varphi_n$ defined by
$$\varphi_n(x)=\left(\left(2-{|x|\over n}\right)\wedge 1\right)^+,\quad \forall\,x\in \R^N.$$ 
It is clear that $u_n$ converges to $u$ in $L$. For the gradients, we have
\[\begin{split}
\int|\nabla(u_n-u)|^2\,dx&\ds\le2\left(\int_{\{n\le|x|\le2n\}}{|u|^2\over n^2}\,dx+\int_{\{n\le|x|\}}|\nabla u|^2\,dx\right)\\
&\le8\int_{\{n\le|x|\}}\left({|u|^2\over|x|^2}+|\nabla u|^2\right)\,dx,
\end{split}\]
which tends to zero as $n\to\infty$.

Let us now consider the case $d=2$. As above, we take $u_n=u\varphi_n$ where $\varphi_n$ is now defined by 
$$\varphi_n(x)=\left({\log\big({2n/|x|}\big)\over\log n}\wedge 1\right)^+,\qquad\forall x\in \R^d.$$
Then $u\varphi_n$ converges to $u$ in $L$ as above and we have
\[\begin{split}
\int|\nabla (u_n-u)|^2\,dx&\le2\int_{\{2\le|x|\le2n\}}\left({|u|^2\over|x|^2\log^2n}
+{\log^2\big({2/|x|}\big)\over\log^2 n}|\nabla u|^2\right)\,dx+\int_{\{2n\le|x|\}}|\nabla u|^2\,dx\\
&\le2\int_{\{2\le|x|\le2n\}}{\log^2|x|\over\log^2n}\left({|u|^2\over|x|^2\log^2|x|}+|\nabla u|^2\right)\,dx+\int_{\{2n\leq |x|\}}|\nabla u|^2\,dx,
\end{split}\]
which tends to zero as $n\to\infty$.\\
In order to show the equivalence of norms, we recall the following Hardy inequalities
\begin{itemize}
\item
If $d=1$, we have
\be\label{hard1}
\int{|u(x)|^2\over|x|^2}\,dx\le4\int|u'(x)|^2\,dx,\qquad\forall u\in C^\infty_c,\ u(0)=0.
\ee
\item
If $d=2$ then, setting $\delta_u=\dist(0,\spt u)$,
\be\label{hard2}
\int{|u(x)|^2\over|x|^2\log^2({|x|/\delta_u})}\,dx\le4\int|\nabla u(x)|^2\,dx,\qquad\forall u\in C^\infty_c\ \hbox{ with }\delta_u>0.
\ee
\item
If $d>2$, then 
\be\label{hard3}
\int{|u(x)|^2\over|x|^2}\,dx\le{4\over(d-2)^2}\int|\nabla u(x)|^2\,dx,\qquad\forall u\in C^\infty_c.
\ee
\end{itemize}
The equivalence with the norm defined by \eqref{norHd>2} for $d\geq 3$ immediately follows from \eqref{hard3}. In the cases $d=1,2$ we take $\varphi\in C^\infty_c$ such that $\varphi=1$ in $B(0,1/2)$, $\varphi=0$ in $\R^d\setminus B(0,1).$ Then, for $u\in H$, we have
$$\|uW\|_{L^2}\leq \|u\varphi W\|_{L^2}+ \|u(1-\varphi)W\|_{L^2}.$$
Since $u\varphi W$ vanishes outside $B(0,1)$, we can apply Poincar\'e's inequality in $H^1_0(B(0,1))$ to deduce
$$ \|u\varphi W\|_{L^2}\leq C\|\nabla (u\varphi W)\|_{L^2}\leq C\big(\|u\|_{L^2(B(0,1))}+\|\nabla u\|_{(L^2)^d}\big),$$
with $C$ independent of $u$. On the other hand, \eqref{hard1} or \eqref{hard2} show the inequality
$$ \|u(1-\varphi)W\|_{L^2}\leq C\big(\|u\|_{L^2(B(0,1))}+\|\nabla u\|_{(L^2)^d}\big).$$
We have thus proved the existence of $C>0$ such that
$$\|u\|_{L}=\|uW\|_{L^2}\leq C\big(\|u\|_{L^2(B(0,1))}+\|\nabla u\|_{(L^2)^d}\big),$$
from which we easily deduce the equivalence of norms.
\end{proof}

\begin{proof}[Proof of Theorem \ref{gacomp}] For every $\O\subset\R^d$ bounded open, it is known the existence of a subsequence of $\mu_n$, still denoted by $\mu_n$, and $\mu_\O\in\M_{cap}$ (in the classical result $\mu_\O$ is only defined in $\O$, but identifying a measure $\mu$ in $\O$ with $\mu\lfloor\O$ we can assume $\mu_\O$ in $\M_{cap}$) such that the sequence of functionals $F_{n,\O}:H^1_0(\O)\to \R$ defined as
$$F_{n,\O}(u)=\int_\O|\nabla u|^2\,dx+\int_\O u^2\,d\mu_n\qquad \forall u\in H^1_0(\O)$$
$\Gamma$-converges in $L^2(\O)$ to $F_\O:H^1_0(\O)\to \R$ defined as
$$F_\O(u)=\int_\O|\nabla u|^2\,dx+\int_\O u^2\,d\mu,\qquad\forall u\in H^1_0(\O).$$
Applying this result to $B(0,k)$ for $k\in \N$ and using a diagonal procedure we can take a subsequence of $n$ such that $F_{n,B(0,k)}$ $\Gamma$-converges for every $k>0$. Since $H^1_0(B(0,k))\subset H^1_0(B(0,m))$ for $k\le m$, we have
$$\int_{B(0,k)}|\nabla u|^2\,dx+\int_{B(0,k)}u^2\,d\mu_{B(0,k)}=\int_{B(0,m)}|\nabla u|^2\,dx+\int_{B(0,m)}u^2\,d\mu_{B(0,m)},$$
for every $u\in H^1_0(B(0,k))$ and therefore
\be\label{ecu1}\mu_{B(0,k)}(E)=\mu_{B(0,m)}(E),\quad \forall\, E\subset B(0,k)\ \hbox{ Borel,}\quad\forall\, k\leq m.\ee
We define $\mu\in \M_{cap}$ by
$$\mu(E)=\lim_{k\to\infty}\mu_{B(0,k)}\big(E\big),$$
where the limit exists because the sequence $\mu_{B(0,k)}\big(E\big)$ is increasing by \eqref{ecu1}.

Let us prove that $\mu_n$ $\gamma$-converges to $\mu$. We consider $u\in H$ and a sequence $u_n\in H$ which converges in $L^2_{loc}$ to $u$. Let us prove
\be\label{ecu2}
\int|\nabla u|^2\,dx+\int|u|^2\,d\mu\le\liminf_{n\to\infty}\left(\int|\nabla u_n|^2\,dx+\int u_n^2\,d\mu_n\right).
\ee
Clearly, we can assume
$$\liminf_{n\to\infty}\left(\int|\nabla u_n|^2\,dx+\int u_n^2\,d\mu_n\right)<\infty,$$
Then, extracting a subsequence if necessary, we have that $u_n$ converges weakly to $u$ in $H^1_{loc}$. For $\varphi\in C^\infty_c$, we have
\[\begin{split}
&\int|\nabla(u\varphi)|^2\,dx+\int|u\varphi|^2\,d\mu\le\liminf_{n\to\infty}
\left(\int|\nabla (u_n\varphi)|^2\,dx+\int|u_n\varphi|^2\,d\mu_n\right)\\
&=\liminf_{n\to\infty}\left(\int|\nabla u_n|^2\varphi^2\,dx+2\int\nabla u_n\cdot\nabla\varphi\,dx+\int|\nabla\varphi|^2\,dx+\int|u_n\varphi|^2\,d\mu_n\right)\\
&\le\liminf_{n\to\infty}\left(\int|\nabla u_n|^2\,dx+\int|u_n|^2\,d\mu_n\right)+
2\int\nabla u\cdot\nabla\varphi\,dx+\int|\nabla\varphi|^2\,dx,
\end{split}\]
which developing the left-hand side shows
$$\int|\nabla u|^2\varphi^2\,dx+\int|u\varphi|^2\,d\mu\le\liminf_{n\to\infty}\left(\int|\nabla u_n|^2\,dx+\int|u_n|^2\,d\mu_n\right).$$
and then \eqref{ecu2} since $\varphi$ is arbitrary.

It remains to prove that for every $u\in H$, there exists a sequence $u_n\in H$ such that
$$\lim_{n\to\infty}\left(\int|\nabla u_n|^2\,dx+\int|u_n|^2\,d\mu_n\right)=\int|\nabla u|^2\,dx+\int |u|^2\,d\mu.$$
If $u$ has compact support the result is immediate from the $\gamma$-convergence of the functionals $F_{n,B(0,k)}$ for every $k$. Therefore it is enough to show that for every $u\in H$ there exists a sequence $u_k$ with compact support which converges to $u$ in $L^2_{loc}\cap L^2_\mu$ and is such that $\nabla u_k\to \nabla u$ in $(L^2)^d$. This is a consequence of the density of $C^\infty_c$ in $H$ proved in Proposition \ref{DenCiH} and a truncation argument.
\end{proof}

\begin{proof}[Proof of Proposition \ref{DePd12}] From \eqref{norHd<=2}, it is enough to show the existence of $C>0$ such that
$$\|u\|^2_{L^2(B(0,1))}\le C\big(\|u\|^2_{L^2_\mu}+\|\nabla u\|^2_{(L^2)^d}\big),\qquad\forall\mu\in\mathcal{U},\ \forall u\in H\cap L^2_\mu.$$
We argue by contradiction. If this inequality does not hold then for every $n\in \N$, there exist $\mu_n\in \mathcal{U}$ and $u_n\in H\cap L^2_{\mu_n}$ such that
$$\|u_n\|_{L^2(B(0,1))}=1,\qquad{1\over n}>\|u_n\|^2_{L^2_{\mu_n}}+\|\nabla u_n\|^2_{(L^2)^d}.$$
In particular, $u_n$ is bounded in $H$ and $\nabla u_n$ converges strongly to zero in $(L^2)^d$. Therefore, up to a subsequence, there exists $u\in H$, constant, such that
$$u_n\rightharpoonup u\hbox{ in }H.$$
Using Theorem \ref{gacomp} and that $\mathcal{U}$ is closed for the $\gamma$-convergence, we can also assume the existence of $\mu\in\mathcal{U}$ such that $\mu_n$ $\gamma$-converges to $\mu$. Thus,
$$u^2\int d\mu=\int u^2d\mu\le\liminf\left(\int|\nabla u_n|^2\,dx+\int u_n^2\,d\mu_n\right)=0.$$
Since the null measure does not belong to $\mathcal{U}$, this shows $u=0$. But this is in contradiction with Rellich-Kondrachov's compactness theorem which implies
$$\|u\|_{L^2(B(0,1))}=\lim_{n\to\infty}\|u_n\|_{L^2(B(0,1))}=1.$$
\end{proof}

\begin{proof}[Proof of Proposition \ref{ConsoEc}] Thanks to Proposition \ref{DenCiH} if $d\ge3$ and Proposition \ref{DePd12} if $d=1,2$, the norm of $u_n$ in $H\cap L^2_{\mu_n}$ is bounded. Thus, for a subsequence of $n$, $u_n$ converges weakly in $H$ to a function $u$. In particular $u_n$ converges weakly in $H^1_{loc}$ and then the classical convergence results for the $\gamma$-convergence show (see e.g. \cite{dmmo87}, \cite{DMGA}) that $u$ satisfies
\be\label{ecu1ab}
\int\nabla u\cdot\nabla v\,dx+\int uv\,d\mu=\langle f,v\rangle,\qquad\forall v\in H^1\cap L^2_\mu\hbox{ with compact support.}
\ee
On the other hand we observe that the inequality
$$\int|\nabla u|^2\,dx+\int|u|^2\,d\mu\le\liminf_{n\to\infty}\left(\int|\nabla u_n|^2\,dx+\int|u_n|^2\,d\mu_n\right),$$
shows that $u$ belongs to $H\cap L^2_\mu$ and then the density of the functions in $H^1\cap L^2_\mu$ with compact support in $H\cap L^2_\mu$ shows that \eqref{ecu1ab} holds for every $v\in H\cap L^2_\mu$. Thus, $u$ is the solution of \eqref{pbun} and then by uniqueness it is not necessary to extract any subsequence.

The convergence of $u_n$ in $W^{1,q}_{loc}$ for $q<2$ is a consequence for example of Proposition 5.4 in \cite{CDGA} where the result is proved for nonlinear systems, the proof is based on the ideas introduced in \cite{BoMu} and \cite{DMMu}.
\end{proof}

\begin{proof}[Proof of Proposition \ref{psisci}] We consider a sequence $(\mu_n)$ in $\M_{cap}$ which $\gamma$-converges to a measure $\mu\in\M_{cap}$. We have to prove that
\be\label{Ec1Lsc}
\Psi(\mu)\le\liminf_{n\to\infty}\Psi(\mu_n).
\ee
Clearly, we can assume that there exists
\be\label{Ec1aLsc}
\lim_{n\to\infty}\Psi(\mu_n)<+\infty.
\ee
We divide the proof in several steps, according to the value of $C_\psi$ defined in \eqref{kpsi}.

\medskip\noindent{\it Case $C_\psi=0$}.
By Remark \ref{Rea0} we can assume that \eqref{supa0} holds. So, for $\eps>0$ we can take $k>0$ such that
\be\label{Ec2Lsc}
\psi(s)\le\eps\qquad\forall\,s\ge k.
\ee
For a fixed $M>0$, we define $\tilde\mu_n\in\M_{cap}$ by setting for every Borel set $E\subset\R^d$
$$\tilde\mu_n(E)=\begin{cases}
\ds\int_E\mu^a_n\wedge k\,dx&\hbox{if }\cp\big(E\setminus B(0,M)\big)=0,\\
+\infty&\hbox{otherwise.}
\end{cases}$$
Possibly extracting a subsequence, we can assume that there exists $h_k\in L^\infty$ such that
\be\label{Ec3Lsc}
\mu^a_n\wedge k\rightharpoonup h_k\quad\hbox{weakly* in }L^\infty,
\ee
and therefore $\tilde\mu_n^a$ $\gamma$-converges to $\tilde\mu$ defined, for every Borel set $E\subset\R^d$, as
$$\tilde\mu(E)=\begin{cases}
\ds\int_E h_k\,dx&\hbox{if }\cp\big(E\setminus B(0,M)\big)=0\\
+\infty&\hbox{otherwise.}
\end{cases}$$
Arguing by convexity thanks to \eqref{Ec3Lsc} and using that \eqref{Ec2Lsc} implies
$$\psi(\mu_n^a\wedge k)\le\psi(\mu_n^a)+\eps,$$
we have
\[\begin{split}
\int_{B(0,M)}\hskip-5pt\psi(h_k)\,dx
&\le\liminf_{n\to\infty}\int_{B(0,M)}\hskip-5pt\psi\big(\mu^a_n\wedge k\big)\,dx\\
&\le\liminf_{n\to\infty}\int_{B(0,M)}\hskip-5pt\big(\psi(\mu_n^a)+\eps\big)\,dx\\
&\le\liminf_{n\to\infty}\Psi(\mu_n)+\eps|B(0,M)|.
\end{split}\]
On the other hand, taking into account that $\mu_n^a\wedge k\le\mu_n$ in $B(0,M)$, we have that $h_k\le\mu$ in $B(0,M)$, which, using that $h_k\in L^\infty$, implies $h_k\le\mu^a$ in $B(0,M)$. Since $\psi$ is decreasing this proves that
$$\int_{B(0,M)}\hskip-5pt\psi(\mu^a)\,dx\le\int_{B(0,M)}\hskip-5pt\psi(h_k)\,dx.$$
Hence,
$$\int_{B(0,M)}\hskip-5pt\psi(\mu^a)\,dx\le\liminf_{n\to\infty}\Psi(\mu_n)+\eps|B(0,M)|,$$
for every $\eps,M>0$. Thus
$$\int\psi(\mu^a)\,dx\le\liminf_{n\to\infty}\Psi(\mu_n),$$
and then \eqref{Ec1Lsc}.

\medskip\noindent{\it Case $0<C_\psi<+\infty$.}
From the definition \eqref{defpsi} of $\Psi$ and \eqref{Ec1aLsc}, we have that the measures $\psi(\mu_n^a)\,dx+C_\psi\mu_n^s$ have uniformly bounded total variation and therefore, for a subsequence, there exists $\nu\in\M$ such that
\be\label{Ec3bLsc}
\psi(\mu_n^a)\,dx+C_\psi\mu_n^s\ \rightharpoonup\nu\ \hbox{ weakly* in }\M.
\ee
Now, we take two bounded open sets $U_2\supset\overline{U_1}$ in $\R^d$ and $\phi\in C^\infty_c(U_2)$ such that $\chi_{U_1}\le\phi\le\chi_{U_2}$. Then, we take a sequence $\phi_n\in H^1_0(U_2)\cap L^2_{\mu_n}(U_2)$ such that
$$0\le\phi_n\le1\hbox{ in }U_2,\qquad\phi_n\rightharpoonup\phi\hbox{ in }H^1_0(U_2),$$
$$\lim_{n\to\infty}\left(\int_{U_2}|\nabla\phi_n|^2\,dx+\int_{U_2}|\phi_n|^2\,d\mu_n\right)=\int_{U_2}|\nabla\phi|^2\,dx+\int_{U_2}|\phi|^2\,d\mu.$$
By Remark \ref{RemRS}, this is equivalent to
\be\label{Ec4Lsc}
\int_{U_2}\nabla\phi_n\cdot\nabla v_n\,dx+\int_{U_2}\phi_nv_n\,d\mu_n\to
\int_{U_2}\nabla\phi\cdot\nabla v\,dx+\int_{U_2}\phi v\,d\mu,
\ee
for every $v\in H^1_0(U_2)\cap L^2_\mu(U_2)$ and every sequence $v_n\in H^1_0(U_2)\cap L^2_{\mu_n}(U_2)$ which converges weakly to $v$ in $H^1_0(U_2)$ and satisfies
$$\limsup_{n\to\infty}\left(\int_{U_2}|\nabla v_n|^2\,dx+\int_{U_2}v_n^2\,d\mu_n\right)<\infty.$$
Now, for a nonnegative $\varphi\in C^1_c(U_1)$, $\varphi\not\equiv 0$ we apply \eqref{desPco} below to deduce
\be\label{Ec5Lsc}
\begin{split}
\int\psi(\mu^a_n)\varphi\,dx+C_\psi\int\varphi\,d\mu^s_n
&\ge\int\psi(\mu^a_n)\phi_n\varphi\,dx+C_\psi\int\phi_n \varphi\,d\mu^s_n\\
&\ge\int\phi_n\varphi\,dx\,\psi\left(\frac{\int\phi_n\varphi\,d\mu_n}{\int\phi_n\varphi\,dx}\right).
\end{split}\ee
By \eqref{Ec4Lsc} we have
$$\int\nabla\phi_n\cdot\nabla\varphi\,dx+\int\phi_n\varphi\,d\mu_n\to\int\nabla\phi\cdot\nabla\varphi\,dx+\int\phi\,\varphi\,d\mu,$$
which combined with $\phi_n$ converging weakly to $\phi$ in $H^1_0(U_2)$ and $\phi=1$ in the support of $\varphi$ shows
$$\int\phi_n\varphi\,d\mu_n\to\int\varphi\,d\mu.$$
Thanks to \eqref{Ec3bLsc} this allows to pass to the limit in $n$ in \eqref{Ec5Lsc} to deduce
$$\frac{\int\varphi\,d\nu}{\int\varphi\,dx}\ge\psi\left(\frac{\int\varphi\,d\mu}{\int\varphi\,dx}\right).$$
for every nonnegative $\varphi\in C^1_c(U_1)$, $\varphi\not\equiv0$. Taking into account that $U_1$ was arbitrary and the derivation measures theorem, we then conclude
\be\label{Ec6Lsc}
\nu^a\ge\psi(\mu^a)\hbox{ in }\R^d.
\ee
On the other hand, for $\eps>0$, Definition \ref{kpsi} of $C_\psi$ proves the existence of $k>0$ such that
\be\label{Ec7Lsc} (C_\psi-\eps)s\le\psi(s)\qquad\forall s\ge k.\ee
Extracting a subsequence we can assume the existence of $h\in L^\infty$ such that
\be\label{Ec8Lsc}
\mu_n^a\wedge k\ \rightharpoonup\ h\qquad\hbox{weakly* in }L^\infty.
\ee
By \eqref{Ec7Lsc}, we have
$$(C_\psi-\eps)\mu_n^a\chi_{\{\mu_n^a>k\}}\le\psi(\mu_n^a)$$
which, by \eqref{Ec3bLsc} and \eqref{Ec8Lsc}, proves
$$(C_\psi-\eps)\mu-(C_\psi-\eps)h\le\nu\qquad\hbox{in }\R^d$$
As a consequence, we obtain that
$$(C_\psi-\eps)\mu^s\le\nu^s\qquad\hbox{in }\R^d$$
and then, by the arbitrariness of $\eps$, that $C_\psi\mu^s\le\nu^s$. This, combined with \eqref{Ec6Lsc}, shows that
$$\psi(\mu^a)+C_\psi\mu^s\le\nu\qquad\hbox{in }\R^d.$$
In particular, for every $\varphi\in C^0_c$, with $0\leq \varphi\leq 1$, we have
\[\begin{split}
\liminf_{n\to\infty}\Psi(\mu_n)
&=\liminf_{n\to\infty}\left(\int\psi(\mu^a_n)\,dx+C_\psi\mu_n^a(\R^d)\right)\\
&\ge\lim_{n\to\infty}\left(\int\psi(\mu^a_n)\varphi\,dx+C_\psi\int\varphi\,d\mu_n^a\right)\\
&=\int\varphi\,d\nu\ge\int\psi(\mu^a)\varphi\,dx+C_\psi\int\varphi\,d\mu_n,
\end{split}\]
which, letting $\varphi$ converge to 1, gives
$$\liminf_{n\to\infty}\Psi(\mu_n)\ge\Psi(\mu).$$
\noindent{\it Case $C_\psi=+\infty$}. In this case \eqref{defpsi} and \eqref{Ec1aLsc} imply that $\mu_n$ is bounded in $L^1$ and equi-integrable in every bounded open set $U$ of $\R^d$ and then the $\gamma$-limit of $\mu_n$ agrees with its weak limit in $L^1(U)$ (see for instance Proposition 2.5 of \cite{bgrv14}). The result is then a consequence of the convexity of $\psi$.
\end{proof}

\begin{lemm}\label{leJg}
We consider $\psi:\R^+\to[0,+\infty]$ convex and lower semicontinuous and a measure $\mu\in\M$. Then, for every $\varphi\in C^0_c$, $\varphi\ge0$ and $\varphi\not\equiv0$, we have
\be\label{desPco}
\int\psi(\mu^a)\varphi\,dx+C_\psi\int\varphi\,d\mu^s\ge\psi\left(\frac{\int\varphi\,d\mu}{\int\varphi\,dx}\right)\int\varphi\,dx\,,
\ee
where $\mu^a$ and $\mu^s$ denote the regular and singular part of $\mu$ and $C_\psi$ is defined in \eqref{decompmu}. Moreover, if in addition $\mu\in\M_{cap}$, then the inequality above also holds for $\varphi\in H^1$, with compact support, nonnegative, and with $\varphi\not\equiv0$.
\end{lemm}

\begin{proof} 
Since $\mu^s$ is singular with respect to the Lebesgue measure, there exists a sequence $C_n$ of measurable sets in $\R^d$ and a sequence of nonnegative measurable functions $h_n$ in $\R^d$ such that denoting
$$i_n={\rm ess\,inf}\,h_n,$$
we have
\be\label{LeJg1}
\lim_{n\to\infty}i_n=\infty,\qquad h_n\chi_{C_n}\rightharpoonup\mu^s\hbox{ weakly* in }\M.
\ee
Now, we consider $\varphi$ in the conditions of the lemma. If $\psi(\mu^a)\varphi$ is not integrable or $C_\psi\varphi$ is not $\mu^s$ integrable, then there is nothing to prove. So, we assume $\psi(\mu^a)\varphi\in L^1$, $C_\psi\varphi\in L^1_{\mu^s}$. Let us show that for every $\ep>0$ we have
\be\label{LeJg2}
\Big(\psi(\mu^a)+(C_\psi +\ep)h_n\chi_{C_n}-\psi(\mu^a+h_n\chi_{C_n})\Big)^-\varphi\to 0\ \hbox{ in }L^1.
\ee
For this purpose we observe that
$$\Big(\psi(\mu^a)+(C_\psi +\ep)h_n\chi_{C_n}-\psi(\mu^a+h_n\chi_{C_n})\Big)\varphi=
\Big(\psi(\mu^a)+(C_\psi +\ep)h_n-\psi(\mu^a+h_n)\Big)\varphi\chi_{C_n},$$
where $\mu^a+h_n\ge i_n$ in $C_n$ with $i_n$ converging to infinity and Definition \ref{kpsi} of $C_\psi$ imply
$$\begin{array}{l}\ds \Big(\psi(\mu^a)+(C_\psi +\ep)h_n-\psi(\mu^a+h_n)\Big)\varphi\chi_{C_n}\\ \noalign{\medskip}\ds\leq\Big(\psi(\mu^a)-(C_\psi +\ep)\mu^a\Big)\varphi\chi_{C_n}\ \hbox{ for }n\hbox{ large enough.}\end{array}$$
Using also that $|C_n\cap\spt(\varphi)|$ tends to zero, which is a consequence of \eqref{LeJg1}, and that $\mu^a$ and $\psi(\mu^a)\in L^1_{loc}$ we deduce that the right-hand side of the previous inequality tends to zero in $L^1$ and then \eqref{LeJg2}.

From \eqref{LeJg2} and Jensen inequality we then have
\[\begin{split}
\int\psi(\mu^a)\varphi\,dx+(C_\psi+\ep)\int\varphi\,d\mu^s
&=\lim_{n\to\infty}\int\big(\psi(\mu^a)+(C_\psi+\ep)\chi_{C_n}h_n\big)\varphi\,dx\\
&\ge\liminf_{n\to\infty}\int\psi(\mu^a+h_n\chi_{C_n})\varphi\,dx\\
&\ge\lim_{n\to\infty}\psi\left(\frac{\int\big(\mu^a+h_n\chi_{C_n}\big)\varphi\,dx}{\int\varphi\,dx}\right)\int\varphi\,dx\\
&=\psi\left(\frac{\int\varphi\,d\mu}{\int\varphi\,dx}\right)\int\varphi\,dx\,.
\end{split}\]
By the arbitrariness of $\ep$, this proves \eqref{desPco}. If we now assume that $\mu$ vanishes on the sets of capacity zero then, using the density of $C^\infty_c(U)$ in $H^1_0(U)$ for every bounded open set $U\subset\R^d$, we get that $\varphi$ can be chosen in $H^1$ with compact support.
\end{proof}

\begin{proof}[Proof of Teorem \ref{ThExsoPo}] Taking into account Proposition \ref{psisci} and that the $\gamma$-limit $\mu$ of a sequence of measures $\mu_n$ satisfying $\mu_n\geq \nu$ also satisfies this restriction, we get that the set $\EE$ of measures satisfying the restrictions
$$\Psi(\mu)\le1,\quad\mu\ge\nu$$
is closed for the $\gamma$-convergence. Therefore, if $\mu_n$ is a minimizing sequence for problem \eqref{optpb}, we get by Theorem \ref{gacomp} that at least for a subsequence, there exists $\mu\in\EE$ such that $\mu_n$ $\gamma$-converges to $\mu$. Observing that if $d=1,2$, condition \eqref{Conpsinu} shows that $\EE$ does not contain the null measure, we can now apply Proposition \ref{ConsoEc} to deduce that the solution $u_n$ of \eqref{pbun} satisfies \eqref{convun} with $u$ the solution of \eqref{pbelimu}. Now, we use
\[\begin{split}
\int j(x,u_n,\nabla u_n)\,dx=&\int\Big(j(x,u_n,\nabla u_n)+g_1|u_n|+g_2|u_n|^2+h_1|\nabla u_n|+h_2|\nabla u_n|^q\Big)\,dx\\
&-\int\Big(g_1|u_n|+g_2|u_n|^2+h_1|\nabla u_n|+h_2|\nabla u_n|^q\Big)\,dx.
\end{split}\]
The first term on the right-hand side of this equality is nonnegative thanks to \eqref{hypoj}. Since \eqref{convun} implies the convergence in measure of $(u_n,\nabla u_n)$, we can then apply Fatou's Lemma to deduce
\[\begin{split}
\liminf_{n\to\infty}&\int\Big(j(x,u_n,\nabla u_n)+g_1|u_n|+g_2|u_n|^2+h_1|\nabla u_n|+h_2|\nabla u_n|^q\Big)\,dx\\
&\ge\int\Big(j(x,u,\nabla u)+g_1|u|+g_2|u|^2+h_1|\nabla u|+h_2|\nabla u|^q\Big)\,dx.
\end{split}\]
Taking into account the convergence in measure of $u_n$ and $\nabla u_n$ and that $|u_n|$ is bounded in $L$ and $|\nabla u_n|$ in $L^2$, we get that $|u_n|$ converges weakly in $H$ to $|u|$ and $|\nabla u_n|$ converges weakly in $L^2$ to $|\nabla u|$. Thus, we have
$$\lim_{n\to\infty} \int\Big(g_1|u_n|+h_1|\nabla u_n|\Big)\,dx=\int\Big(g_1|u|+h_1|\nabla u|\Big)\,dx.$$
On the other hand, the strong convergences of $u_n$ in $L^2_{loc}$ and $\nabla u_n$ in $(L^q_{loc})^d$ imply
$$\lim_{n\to\infty}\int_{B(0,R)}\Big(g_2|u_n|^2+h_2|\nabla u_n|^q\Big)\,dx=\int_{B(0,R)}\Big(g_2|u|^2+h_2|\nabla u|^q\Big)\,dx,$$
for every $R>0$, while estimate
\[\begin{split}
\int_{\R^d\setminus B(0,R)}\Big(g_2|u_n|^2+h_2|\nabla u_n|^q\Big)\,dx\le&\left\|{g_2\over W^2}\right\|_{L^\infty(\R^d\setminus B(0,R))}\|u_n\|_{L}^2\\
&+\|h_2\|_{L^{2/(2-q)}(\R^d\setminus B(0,R))}\|\nabla u_n\|^q_{(L^2)^d}
\end{split}\]
combined with assumption \eqref{hypoab}, $u_n$ bounded in $L$ and $|\nabla u_n|$ bounded in $L^2$ show 
$$\lim_{R\to\infty}\limsup_{n\to\infty}\int_{\R^d\setminus B(0,R)}\Big(g_2|u_n|^2+h_2|\nabla u_n|^q\Big)\,dx=0,$$
and therefore
$$\lim_{n\to\infty}\int\Big(g_2|u_n|^2+h_2|\nabla u_n|^q\Big)\,dx=\int\Big(g_2|u|^2+h_2|\nabla u|^q\Big)\,dx.$$
We have then proved
$$\liminf_{n\to\infty}\int j(x,u_n,\nabla u_n)\,dx\ge\int j(x,u,\nabla u)\,dx$$
and thus that $\mu$ is a solution of \eqref{optpb}.
\end{proof}

\section{Some necessary conditions of optimality}\label{snece}

In the previous sections, using the $\gamma$-convergence theory we have studied the existence of solution for problem \eqref{optpb}. Here let us show how assuming some derivability conditions for the functions $j$ and $\psi$ defining the
cost function and the volume restriction respectively, we can obtain some optimality conditions for \eqref{optpb}. Namely, for $j=j(x,s,\xi)$ let us assume the existence of $\partial_sj(x,s,\xi)$, $\partial_\xi j(x,s,\xi)$ where these functions are continuous in $(s,\xi)$, measurable in $x$ and satisfy the growth condition
\be\label{grcodej} {\big|\partial_sj(x,s,\xi)\big|\over W}+ |\partial_\xi j(x,s,\xi)\big|\leq k+M\big(W|s|+|\xi|\big)\quad \forall\, (s,\xi)\in \R\times\R^d,\ \hbox{a.e. }x\in \R^d,\ee
with $k\in L^2$, $M>0$.

For the function $\psi$ we assume that it is finite in $(0,\infty)$ and continuous and derivable in $[0,\infty)$. Here, if $\psi(0)=\infty$, we define $\psi'(0)=-\infty$. Indeed, we observe that in this case the continuity of $\psi$ in $0$ and its convexity imply
$$\lim_{s\to 0^+}\psi=+\infty,\qquad \lim_{s\to 0^+}\psi'(s)=-\infty.$$
With these conditions, the following result holds.
\begin{theo} \label{thopco} In the assumptions of Theorem \ref{ThExsoPo}, we assume that $j$ and $\psi$ in problem \eqref{optpb} satisfy the conditions stated above and that there exists a measure $\tilde \mu\in \M_{cap}$ such that 
\be\label{hicoo}\tilde \mu\geq \nu,\qquad \Psi(\tilde \mu)<1,\ee
Then, if $\mu\in \M_{cap}$ is a solution of problem \eqref{optpb}, $u$ is the corresponding state function and $p$, the adjoint state, is defined as the solution of
\be\label{defAdSt} -\Delta p+p\mu=\partial_sj(x,u,\nabla u)-{\rm div}\, \partial_\xi j(x,u,\nabla u)\ \hbox{ in }\R^d,\quad p\in H\cap L^2_\mu,\ee
we have the existence of $\lambda\geq 0$ such that
\be\label{Con1opco} \lambda\big(\Psi(\mu)-1)=0,\ee
\be\label{Con2opco} \lambda\psi'(\mu^a)\geq up\ \hbox{ a.e. in }\R^d,\ee
\be\label{Con3opco} \lambda\psi'(\mu^a)=up\ \hbox{ a.e. in }\big\{\mu^a>\nu^a\big\}.\ee
Moreover, if $C_\psi$ defined by \eqref{kpsi} is finite, we have:
\be\label{Con4opco} \lambda C_\psi\geq up\ \hbox{ q.e. in }\R^d,\qquad(\mu^s-\nu^s)\big(\{\lambda C_\psi>up\}\big)=0.\ee
\end{theo}

\begin{rema} Observe that $p$ is well defined thanks to $u\in H$ and \eqref{grcodej} which imply that the right-hand side in \eqref{Con1opco} is in $H'$.
\end{rema}

\begin{rema} Assume that in Theorem \ref{thopco} the constant $\lambda$ is positive and $\psi'$ is strictly increasing, then conditions \eqref{Con2opco}, \eqref{Con3opco} provide
$$\mu^a=\left\{\begin{array}{ll} \nu^a &\hbox{ if }\lambda\psi'(\nu^a)>up\\ \ecart\ds
(\psi')^{-1}\Big({up\over \lambda}\Big) &\hbox{ if }\lambda\psi'(\nu^a)\leq up.\end{array}\right.$$
\end{rema}

As a consequence of Theorem \ref{thopco} we can now obtain the following result providing sufficient conditions to have a solution $\mu$ of \eqref{optpb} such that $(\mu-\nu)^s=0$ or $\mu=\nu+\infty_K$ with $K$ quasi-closed.

\begin{theo}\label{SoespCC}
In the assumptions of Theorem \ref{thopco} and assuming $C_\psi=0$ and $j=j(x,s,\xi)$ concave in $(s,\xi)\in\R\times\R^d$, there exists a solution $\hat\mu$ of \eqref{optpb} such that $(\hat\mu-\nu)^s=0$. Moreover, if there exists a solution $\mu$ of \eqref{optpb} such that the constant $\lambda$ in Theorem \ref{thopco} vanishes, then we can take $\hat\mu$ as $\hat\mu=\nu+\infty_K$ with $K$ a quasi-closed set of $\R^d$.
\end{theo}

\begin{rema} The concavity condition for $j$ is not very usual in optimization where it is most frequent to deal with $j$ convexe in $s$ and $\xi$ but then it is simple to check that the above result does not hold. Just consider $f\in L$ with $f>0$ in $\R^d$, $\mu_0\in \M_{cap}$ such that $\mu_0\geq \nu$, $\Psi(\mu_0)\leq 1$ and $u_0$ the solution of 
$$-\Delta u_0+\mu_0 u_0=f\ \hbox{ in }\R^d,\quad u_0\in H\cap L^2_{\mu_0}.$$ 
Then, similarly to Lemma 3.3 in \cite{DMGA} it is possible to check that $\mu_0$ is univocally determined by $u_0$. That is, if $\mu\in \M_{cap}$ is such that $u_0$ satisfies
$$-\Delta u_0+\mu u_0=f\ \hbox{ in }\R^d,\quad u_0\in H\cap L^2_\mu,$$ 
then $\mu=\mu_0$. Therefore, if we consider problem \eqref{optpb} with $j=j(x,s)$ the convex function in $s$ defined by
$$j(x,s)=|s-u_0(x)|^2W(x)^2\qquad\forall s\in\R,\hbox{ a.e. in }\R^d$$
we get that $\mu_0$ is the unique solution of \eqref{optpb}. Thus, we cannot ask $\mu_0$ to satisfy any additional assumption to those given by the restrictions in \eqref{optpb}.\end{rema}
\begin{rema} An interesting application of Theorem \ref{SoespCC} corresponds to $j$ linear in $(s,\xi)$. In particular, we can take $f=g-{\rm div}\,G$ with $W^{-1}g\in L^2$, $G\in (L^2)^d$ and
$$j(x,s,\xi)=g(x)s+G(x)\cdot\xi\quad\hbox{or}\quad j(x,s,\xi)=-g(x)s-G(x)\cdot\xi,$$
which respectively correspond to the minimization and the maximization of the energy, i.e. to problems
$$\min\left\{\int|\nabla u|^2\,dx+\int|u|^2\,d\mu\ :\ -\Delta u+\mu u=f,\ \Psi(\mu)\le1,\ \mu\ge\nu\right\}$$
and
$$\max\left\{\int|\nabla u|^2\,dx+\int|u|^2\,d\mu\ :\ -\Delta u+\mu u=f,\ \Psi(\mu)\le1,\ \mu\ge\nu\right\}.$$\end{rema}
\begin{rema} The assumption $\lambda=0$ in Theorem \ref{SoespCC} always holds in the case $\Psi=0$, i.e. where there is not a ``volume'' restriction.
\end{rema}

We finish this section with the following result relative to the case $f\ge0$.

\begin{prop}\label{vradon}
If $\mu\in\M_{cap}$ is a Radon measure and $f\ge0$ with $f\not\equiv0$, then the solution $u$ of the PDE
\be\label{ecuesPp} -\Delta u+\mu u=f\ \hbox{ in }\R^d,\quad u\in H\cap L^2_\mu.\ee
satisfies $u>0$ q.e. in $\R^d$. In addition, when $C_\psi>0$, $f\ge0$, $f\not\equiv 0$, the set
$$\mathcal{U}=\big\{u\in H\ :\ \exists\mu\in \M_{cap}\ \hbox{ with }u\in L^2_\mu,\ -\Delta u+\mu u=f\ \hbox{ in }\R^d,\ \mu\geq \nu,\ \Psi(\mu)\leq 1\big\},$$
is a convex set. 
\end{prop}

\begin{rema}
The assumption that $\mu$ is a Radon measure is essential in Proposition \ref{vradon}; indeed, take in $\R^d$ the measure
$$\mu(E)=\int_{E\cap B}\frac{1}{|x_1|}\,dx$$
where $B$ is the unit ball of $\R^d$ centered at the origin. The measure $\mu$ belongs to $\M_{cap}$ but it is not a Radon measure; in addition, every $u\in H\cap L^2_\mu$ must vanish on the set $\{x\in B\ :\ x_1=0\}$ which is not of capacity zero.
\end{rema}

\begin{rema}
Writing \eqref{optpb} as
$$\min_{u\in \mathcal{U}}\int j(x,u,\nabla u)\,dx,$$
and assuming $j$ strictly convex in $(s,\xi)$, we deduce from Theorem \ref{ThExsoPo} the existence and uniqueness of an optimal state function $\hat u$ and then of an optimal measure 
$$\hat \mu={1\over \hat u}\big(f+\Delta \hat u\big).$$
\end{rema}

\section{Proofs of the results of section \ref{snece}}\label{sprooCO}

\begin{proof}[Proof of Theorem \ref{thopco}] We define 
$$\Theta=\left\{\vartheta\in\M_{cap}\ :\ \vartheta\geq\nu,\ \exists\, r_1\in L^\infty,\ r_2\in L^\infty_\mu\ \hbox{ with }\vartheta=r_1W^2+r_2\mu\right\}.$$
For $\vartheta\in\Theta$, such that $\Psi(\vartheta)\le1$ and $\ep\in[0,1)$ we take $\mu_\ep=(1-\ep)\mu+\ep\vartheta\in \M_{cap}$ and $u_\ep$ as the solution of the corresponding state problem
$$-\Delta u_\ep+\mu_\ep u_\ep=f\ \hbox{ in }\R^d,\qquad u_\ep\in H\cap L^2_{\mu_\ep}.$$
By definition of $\Theta$, we know that
there exist $r_1\in L^\infty$, $r_2\in L^\infty_\mu$ such that $\vartheta=r_1W^2+r_2\mu$. Then, it is simple to check that $H\cap L^2_{\mu_\ep}=H\cap L^2_\mu$ and therefore the equation for $u_\ep$ can be written as 
$$ -\Delta u_\ep+\mu u_\ep+\ep(\vartheta-\mu)u_\ep=f\ \hbox{ in }\R^d,\qquad u_\ep\in H\cap L^2_\mu.$$
This allows us to prove that the function $\ep\in[0,1)\mapsto u_\ep\in H\cap L^2_\mu$ is derivable on the right at cero. Namely, we have
\be\label{ecuOC1} \lim_{\ep\to 0^+}{u_\ep-u\over \ep}=u'\ \hbox{ in }H\cap L^2_\mu,\ee
with $u'$ the solution of
\be\label{ecuOC1b} -\Delta u'+\mu u'+(\vartheta-\mu)u=0\ \hbox{ in }\R^d,\qquad u'\in H\cap L^2_{\mu}.\ee
Now, we use that by convexity $\mu_\ep$ also satisfies the restrictions $\Psi(\mu^\ep)\leq 1$, $\mu_\ep\geq \nu$ and thus
$$\int j(x,u_\ep,\nabla u_\ep)\,dx\ge\int j(x,u,\nabla u)\,dx,\quad\forall\,\ep\in [0,1).$$
Deriving on the right on $\ep=0$ thanks to assumptions \eqref{grcodej} and \eqref{ecuOC1}, we deduce
$$\int \big(\partial_s j(x,u,\nabla u)u'+\partial_\xi j(x,u,\nabla u)\cdot \nabla u'\big)\,dx\ge0,$$
which using $u'$ as test function in \eqref{defAdSt} and then $p$ as test function in \eqref{ecuOC1b} can be written as
\[\begin{split}
\int\big(\partial_s j(x,u,\nabla u)u'+\partial_\xi j(x,u,\nabla u)\cdot \nabla u'\big)\,dx\\
=\int\nabla p\cdot\nabla u'+\int pu'\,d\mu=\int up(d\mu-d\vartheta).
\end{split}\]
We have thus proved
$$\int up\,d\vartheta\le \int up\, d\mu,\qquad\forall\vartheta\in\Theta\hbox{ with }\Psi(\vartheta)\le1.$$
Since $\Theta$ is a convex set and the function $\vartheta\to\int up\,d\vartheta$ is linear (and then convex) we can apply Kuhn-Tucker's theorem to deduce the existence of $\lambda_0,\lambda\geq 0$ not simultaneously zero such that \eqref{Con1opco} is satisfied and
\be\label{teomu1}
-\lambda_0\int up\, d\vartheta+\lambda \Psi(\vartheta)\geq -\lambda_0\int up\, d\mu+\lambda \Psi(\mu),\quad \forall\, \vartheta\in \Theta.
\ee
Let us prove that $\lambda_0$ cannot be zero. We reason by contradiction, if $\lambda_0=0$ then $\lambda\neq0$ and then
$$\Psi(\vartheta)\geq 1=\Psi(\mu)\quad \forall\, \vartheta\in \Theta.$$
Assume $\vartheta\in \M_{cap}$ such that $\vartheta^a\geq \nu^a$, $\vartheta^s+\vartheta^\infty=\nu^s+\nu^\infty$ and define
$$Z_n=\Big\{x\in\R^d: {\vartheta^a(x)\over W(x)^2}\leq n\Big\},\quad 
\vartheta_n=\vartheta^a\chi_{Z_n}+\mu^a\chi_{\R^d\setminus Z_n}+\nu^s+\nu^\infty.$$
Then, $\vartheta_n$ is in $\Theta$ and therefore
$$1\leq \Psi(\vartheta_n)=\int_{Z_n}\psi (\vartheta^a)\,dx+\int_{\R^d\setminus Z_n}\psi(\mu^a)+C_\psi\nu^s(\R^d).$$
Using that $Z_n$ increases to $\R^d$ and that $\psi(\mu^a)$ is in $L^1$ we can use the monotone convergence theorem in the first integral and the Legesgue dominated convergence theorem in the second one to pass to the limit in $n$, obtaining
$$1\leq \int\psi(\vartheta^a)+C_\psi\nu^s(\R^d),\quad\forall\, \theta\in \M_{cap}\hbox{ with }\vartheta^a\geq \nu^a,\ \vartheta^s+\vartheta^\infty=\nu^s+\nu^\infty,$$
but this is a contradiction with \eqref{hicoo}, which taking $\vartheta=\tilde\mu^a+\nu^s+\nu^\infty$ provides
$$1>\Psi(\tilde\mu)=\int\psi(\tilde\mu^a)+C_\psi\tilde \mu^s\geq \int\psi(\tilde\mu^a)+C_\psi\nu^s=\Psi(\vartheta)\geq 1.$$
This contradiction shows that $\lambda_0$ is not zero. Dividing by $\lambda_0$, we can then assume $\lambda_0=1$ in \eqref{teomu1}. Taking into account the definition \eqref{defpsi} of $\Psi$ and the convexity and derivability assumptions on $\psi$ we get that \eqref{teomu1} is equivalent to
$$-\int up\,d\vartheta+\lambda\int \psi'(\mu^a)\,d\vartheta^a+\lambda C_\psi\vartheta^s(\R^d)\geq-\int up\,d\mu+\lambda\int \psi'(\mu^a)\,d\mu^a+\lambda C_\psi\mu^s(\R^d),$$
for every $\vartheta\in \Theta$. An approximation argument similar to the one used above to prove $\lambda_0\not=0$ allows us to take $\vartheta$ such that $\vartheta^a\geq \nu^a$, $\vartheta^s$ absolutely continuous with respect to $\mu^s$, $\vartheta^s\geq \nu^s$ and $\vartheta^\infty=\mu^\infty$. Then \eqref{teomu1} provides
$$\int\big(\lambda\psi'(\mu^a)-up\big)\,d\vartheta^a\ge\int\big(\lambda\psi'(\mu^a)-up\big)\,d\mu^a,\quad\vartheta^a\ge\nu^a,$$
$$\int \big(\lambda C_\psi-up\big)d\vartheta^s\geq \int \big(\lambda C_\psi-up\big)d\mu^s,\quad \vartheta^s\geq \nu^s,$$
with $\vartheta^a$ absolutely continuous with respect to the Lebesgue measure and $\vartheta^s$ absolutely continuous with respect to $\mu^s$. These two conditions are equivalent to \eqref{Con2opco}, \eqref{Con3opco} and \eqref{Con4opco}.
\end{proof}

\begin{proof}[Proof of Theorem \ref{SoespCC}] Let $\tilde \mu$ be a solution of problem \eqref{optpb} and $u$ be the corresponding state function. We observe that defining $\mu$ as $\mu=\tilde \mu+\infty_{\{u=0\}}$, we get that $u$ is also a solution of
$$-\Delta u+\mu u=f\ \hbox{ in }\R^d,\quad u\in H\cap L^2_{\mu}.$$
Moreover, $\mu\geq\tilde\mu\geq \nu$ in $\R^d$ and since by Remark \ref{Rea0}, $C_\psi=0$ implies $\psi$ decreasing, we also have $\Psi(\mu)\leq \Psi(\tilde\mu)\leq 1.$ Therefore, $\mu$ is a solution of \eqref{optpb} with the same state function $u$, which is strictly positive $(\mu^a+\mu^s)$-a.e. in $\R^d$.

Now, we take $p$ as the solution of \eqref{defAdSt} and we observe that \eqref{Con4opco}, $C_\psi=0$ and $u\not=0$ $\mu^s$-a.e. in $\R^d$ imply 
\be\label{CaCv1} up\le0\ \hbox{ q.e. in }\ \R^d, \qquad
p=0\ \ (\mu^s-\nu^s)\hbox{-a.e. in }\R^d.\ee
Moreover, if $\lambda=0$, then condition \eqref{Con2opco} also gives
\be\label{CaCv2} p=0\ \ (\mu^a-\nu^a)\hbox{-a.e. in }\R^d.\ee
We define $\hat \mu$ as
$$\hat\mu=\mu+\infty_{\{p=0\}}.$$
Then, by \eqref{CaCv1} and \eqref{CaCv2} we have that $\hat\mu$ is such that $(\hat\mu-\nu)^s=0$ and for $\lambda=0$, $\hat \mu=\hat\nu+\infty_{\{p=0\}}$, where the set $\{p=0\}$ is quasi-closed.

Let us prove that $\hat\mu$ is also a solution of \eqref{optpb}. First we use that $\psi$ decreasing and $\hat\mu\ge\mu$ proves $\Psi(\hat\mu)\le\Psi(\mu)\le1$.

Now, we define $\hat u$ as the solution of
\be\label{CaCv3} -\Delta \hat u+\hat \mu \hat u=f\ \hbox{ in }\R^d,\quad \hat u\in H\cap L^2_{\hat\mu}.\ee
Since $j(x,s,\xi)$ is concave in $(s,\xi)$ we have
$$j(x,u,\nabla u)\ge j(x,\hat u,\nabla\hat u)+\partial_sj(x,u,\nabla u)(u-\hat u)+\partial_\xi j(x,u,\nabla u)\cdot\nabla(u-\hat u)\ \hbox{ a.e. in }\R^d,$$
and then
\be\label{CaCv4}
\begin{split}
\int j(x,u,\nabla u)\,dx\ge\int j(x,\hat u,\nabla\hat u)\,dx\\
+\int\Big(\partial_sj(x,u,\nabla u)(u-\hat u)+\partial_\xi j(x,u,\nabla u)\cdot\nabla(u-\hat u)\Big)\,dx.
\end{split}\ee
In the last integral we observe that $\hat\mu\geq \mu$ and $\hat u\in L^2_{\hat\mu}$ imply $\hat u\in L^2_\mu$. Then, we can take $u-\hat u$ as test function in \eqref{defAdSt} to get
\[\begin{split}
\int\hskip-3pt\Big(\partial_sj(x,u,\nabla u)(u-\hat u)+\partial_\xi j(x,u,\nabla u)\cdot\nabla (u-\hat u)\Big)\,dx\\
=\int\nabla p\cdot\nabla(u-\hat u)\,dx+\int(u-\hat u)p\,d\mu,
\end{split}\]
but taking $p$ as test function in the equation for $u$ we have
$$\int\nabla p\cdot\nabla u\,dx+\int pu\,d\mu=\langle f,p\rangle_{H',H}.$$
The definition of $\hat\mu$ shows that $p$ belongs to $L^2_{\hat\mu}$, thus we can also take $p$ as test function in \eqref{CaCv3} to get
$$\int\nabla p\cdot\nabla \hat u\,dx+\int p\hat u\,d\hat \mu=\langle f,p\rangle_{H',H},$$
where
$$\int p\hat u\,d\hat \mu=\int p\hat u\,d\infty_{\{p=0\}}+\int p\hat u\,d\mu=\int p\hat u\,d\mu.$$
This proves
$$\int\hskip-3pt\Big(\partial_sj(x,u,\nabla u)(u-\hat u)+\partial_\xi j(x,u,\nabla u)\cdot\nabla (u-\hat u)\Big)\,dx=0,$$
which substituted in \eqref{CaCv4} and using that $u$ is the state function associated to $\mu$ solution of \eqref{optpb} shows that $\hat\mu$ is also a solution of \eqref{optpb}.
\end{proof}

\begin{proof} [Proof of Proposition \ref{vradon}] In order to prove that $u>0$ q.e. in $\R^d$, we first use $u^-$ as test function in (\ref{ecuesPp}) which provides $u\geq 0$ q.e. in $\R^d$. Then, for $\varphi\in C^\infty_c$ and $\ep>0$, we take $(\ep-u)^+\varphi^2$ as test function in (\ref{ecuesPp}). Denoting
$$A_\ep=\big\{x\in \R^d:\quad 0<u<\ep\big\},$$
we get
$$-\int_{A_\ep}|\nabla u|^2\varphi^2\,dx+2\int_{A_\ep}\nabla u\cdot\nabla\varphi\,(\ep-u)\varphi\,dx+\int_{A_\ep}u(\ep-u)\varphi^2\,d\mu\ge0,$$
which using Young's inequality provides
$${1\over 2}\int_{A_\ep} |\nabla u|^2\varphi^2\,dx\le2\ep^2\int_{A_\ep}|\nabla\varphi|^2\,dx+\ep^2\int_{A_\ep}\varphi^2\,d\mu.$$
Passing to the limit as $\ep\to0$, thanks to the fact that $\chi_{A_\ep}$ tends to zero q.e. in $\R^d$, we then deduce
\be\label{limep}
\lim_{\ep\to 0}\int_{A_\ep}|\nabla u|^2\varphi^2\,dx=0.
\ee
Now, we introduce 
$$z_\ep=\Big({u\over\ep}\wedge 1\Big)\varphi\in H,$$
and we observe that
$$z_\ep\to\varphi\chi_{\{u>0\}}\quad\hbox{ in }L,$$
while by (\ref{limep})
$$\nabla z_\ep={1\over \ep}\nabla u\,\varphi\chi_{A_\ep}+\Big({u\over\ep}\wedge 1\Big)\nabla\varphi\to \nabla\varphi\chi_{\{u>0\}}\quad\hbox{ in }(L^2)^d.$$
So,
$$\varphi\chi_{\{u>0\}}\in H,\quad \nabla\big(\varphi\chi_{\{u>0\}}\big)=\nabla\varphi\chi_{\{u>0\}} \hbox{ in }\R^d,\qquad\forall\,\varphi\in C^\infty_c,\ \varphi\ge0,$$
but on the other hand, since $\varphi$ is smooth we have
$$\nabla \big(\varphi\chi_{\{u>0\}}\big)=\nabla \varphi\chi_{\{u>0\}}+\varphi\nabla \chi_{\{u>0\}},$$
and thus
$$\varphi\nabla \chi_{\{u>0\}}=0,\qquad\forall\,\varphi\in C^\infty_c.$$
This proves that $\chi_{\{u>0\}}$ is a constant function in $\R^d$, which can only take the values one or zero. If it is zero then $u$ is the null function but then \eqref{ecuesPp} implies $f=0$ in $\R^d$ contrary to the assumptions of the Proposition. So $\chi_{\{u>0\}}$ is the constant function equals to one in $\R^d$. This proves that the set $\{u=0\}$ has null measure. Let us now see that in fact it has zero capacity. For this purpose we take a compact set $K\subset\{u=0\}$ and $\varphi\in C^\infty_c$, such that $\varphi\ge\chi_K$. Then defining for $\ep>0$ the function $\tilde z_\ep$ as
$$\tilde z_\ep=\Big(1-{u\over\ep}\Big)^+\varphi\in H^1,$$
we have that $z_\ep\ge\chi_K$ q.e. in $\R^d$ and then
$$\begin{array}{ll}\displaystyle
\cp(K)&\displaystyle\le\int|\nabla z_\ep|^2\,dx+\int|z_\ep|^2\,dx\\
\noalign{\medskip}&\displaystyle\le{2\over\ep^2}\int_{A_\ep}|\nabla u|^2\varphi^2\,dx+2\int_{\{0\le u< \ep\}}|\nabla\varphi|^2\,dx+\int_{\{0\leq u< \ep\}}|\varphi|^2\,dx.
\end{array}$$
Letting $\ep$ tend to zero and using (\ref{limep}) and that $\big|\{u=0\}\big|=0$ we deduce that the right-hand side tends to zero. This proves that $\cp(K)=0$ for every compact set $K\subset\{u=0\}$ and then that $\cp(\{u=0\})=0$.

For the convexity of the set $\mathcal{U}$ observe that if $u_1$, $u_2$ belong to $\mathcal{U}$ with respectively associated measures $\mu_1$, $\mu_2$, then for every $\lambda\in (0,1)$, the function $u_\lambda=\lambda u_1+(1-\lambda)u_2$
satisfies 
$$-\Delta u_\lambda+\mu_\lambda u=f\ \hbox{ in }\R^d,\quad u_\lambda\in H\cap L^2_{\mu_\lambda},$$
with
$$\mu_\lambda :=\lambda{u_1\over u_\lambda}\mu_1+(1-\lambda){u_2\over u_\lambda}\mu_2.$$
Moreover $\mu_\lambda$ also satisfies the restrictions $\mu_\lambda\ge\nu$ and $\Psi(\mu_\lambda)\le1$.
\end{proof}

\section{Some numerical simulations}\label{snume} 

In this section, we present some numerical experiments that illustrate some of the qualitative properties of the optimal solutions and the fact that the problem is posed in the whole space $\R^d$. Our numerical simulations are made in the case $d=2$.

For our experiments, having in mind the assumptions on $j(x,s,\xi)$ in the existence Theorem \ref{ThExsoPo} and Theorem \ref{SoespCC} of optimality conditions, we consider the linear case $j(x,s,\xi)=gu$ with $W^{-1}g\in L^2$. With respect to the volume constraint, we consider two different functions:
$$\begin{cases}
\psi(s)=\psi_1(s)=\frac{1}{m}e^{-\alpha s}\ \hbox{ for some }\alpha,m>0\\
\psi(s)=\psi_2(s)=s^2.
\end{cases}$$
For the first choice of $\psi$ one has $C_\psi=\psi(\infty)=0$, while for the second one $C_\psi=\psi(\infty)=+\infty$.
In the first case the volume constraint $\Psi(\mu)\le 1$ reduces to 
$$\int \frac{1}{m}e^{-\alpha\mu^a}\,dx\le1,$$
while in the second one it gives
$$\mu^s=\mu^\infty=0,\quad \int |\mu^a|^2\,dx\le1.$$
Then our goal si to solve numerically problems in the form:

\be\label{optpbnum}
\min\left\{\int g u\,dx\ :\ \int\psi(\mu^a)\,dx\le1,\ \mu\ge\nu\right\}.
\ee
for $u$ the solution of the state equation
\be\label{eq:statenum}
-\Delta u+\mu u=f\ \hbox { in }\R^d,\quad u\in H\cap L^2_\mu.
\ee
For the first case, since $\psi_1(0)>0$, we can assume $\nu\equiv0$, dropping the constraint $\mu\ge\nu$. For the case of $\psi_2$, since we are in the case $d=2$, we have to impose the constraint $\mu\ge\nu$ with $\nu$ different to the null measure. Taking into account Theorem \ref{SoespCC}, in the first case, we can just look for an optimal control $\mu$ of the form $\mu=\mu^a+\mu^\infty$. In fact, numerically let us search for $\mu:\R^d\to[0,+\infty]$ a Lebesgue measurable function.

For the implementation of our algorithm the main required data are the initialization $\mu_0$, the associated routines to the cost and volume function and the associated routines to the gradient of the cost and volume function using the adjoint state. The admissible measures, or equivalently potentials $V$, take values in $[0,+\infty]$. From the numerical point of view it is advisable to constrain $\mu$ to take values on a bounded interval $[0,\mu_{max} ]$, with $\mu_{max}$ large enough. Analogously, in order to solve numerically the extremal problem \eqref{optpbnum} we replace $\R^2$ by $D=(-M, M)\times(-M, M)$ for different values of the constant $M$. This will allow us to study the behavior of the solution of the optimization problem as $M$ goes to $+\infty$.

In order to solve numerically the problem \eqref{optpbnum}, we will use a gradient descent algorithm. Then it is easy to check that the derivative of the cost functional at $\mu$ in the direction $\mu_1$ is given by:
$$\delta I(\mu)[\mu_1]=\int\mu_1(x)u(x)p(x)\,dx,$$
where $p$ is the unique solution of the adjoint system \eqref{defAdSt}.

\subsection{Numerical examples}

For our numerical experiments we decided to use the free software FreeFEM++ v 3.56 (see {\tt http://www.freefem.org/}, see \cite{FF}), complemented with the library NLopt (see {\tt http://ab-initio.mit.edu/wiki/index.php/NLopt}) using the Method of Moving Asymptotes as the optimizing routing (see \cite{svanberg87}). This technique is a gradient method based on a spatial type of convex approximation where in each iteration a strictly convex approximation subproblem is generated and solved. We use this algorithm after to get satisfactory results for a similar problem in the case for bounded domains. This algorithm was previously used in \cite{bmv18} and in \cite{bebuve} in the case when the state equation is posed in a bounded domain.

The minimization problem \eqref{optpbnum} is posed for capacitary measures defined on the whole $\R^2$. Having in mind the definition of the space $L$ in \eqref{defiL}, and the fact that $W^{-1}g\in L^2$ we take for our experiments the approximations of the constant function $g=1$ given by:
$$g(x,y)=\frac{1}{1+\ep(x^2+y^2)^3},$$
with $\ep>0$ a small parameter, for instance, we have considered $\ep=10^{-10}$.

We use a $P_2$-Lagrange finite element approximations for $u$ and $p$, solutions of the state and co-state equations \eqref{eq:statenum} and \eqref{defAdSt}, respectively, and $P_0$ Lagrange finite element approximation for the capacitary measure $\mu$ with $\mu_{max}=15000$. We consider a regular mesh where the number of elements increases with $M$. For instance, in the case of $M=5$ we take a mesh with $2000$ elements while for $M=20$ we consider a mesh with $8000$ elements. We present different results concerning the two examples of functions $\psi$ given above, and different choices of the right-hand side $f$ in the state equation.

{\bf Case 1:} $\psi(s)=\frac{1}{m}e^{-\alpha s}$. We remind a result in \cite{bgrv14}, where the choice of $\psi(s)=e^{-\alpha s}$ is proposed to approximate shape optimization problems with Dirichlet condition on the free boundary. As $\alpha\to0$, the problem approximates the shape optimization problem
$$
\min\left\{\int_\O g u\,dx\ :\ -\Delta u=f\hbox{ in }\O,\ u\in H,\ u=0\hbox{ q.e. in }\R^d\setminus\O,\ |\O|\le m\right\}.
$$
We fix, in our simulation, the value of the parameter $\alpha$ as $\alpha=3\cdot10^{-4}$. 

\begin{figure}[!t]
\begin{minipage}[!t]{14cm}
\centering
\includegraphics[width=6.8cm]{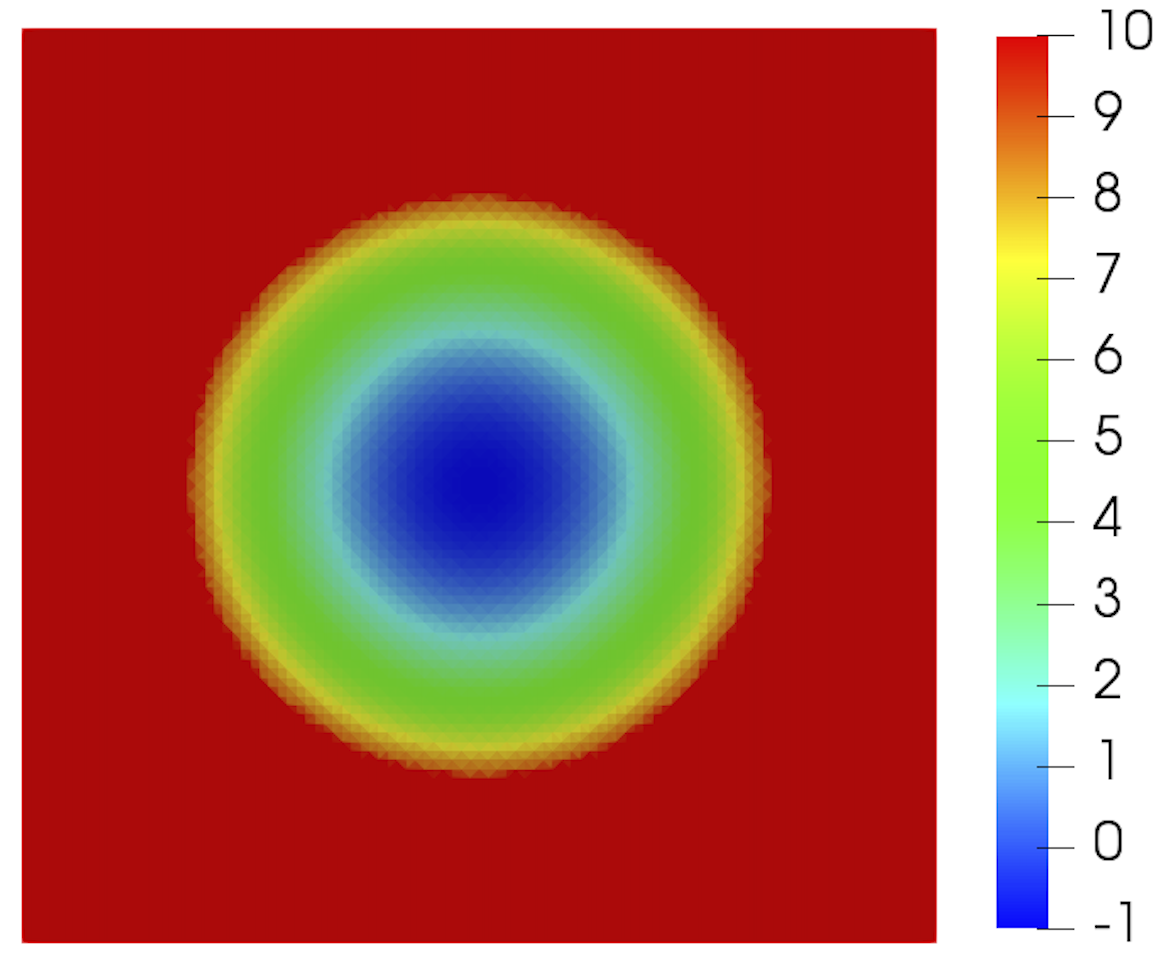}
\end{minipage}\vskip-0.2cm
\caption{Right-hand side function $f_1$ in $D=(-5,5)\times(-5,5)$ (right).}\label{F1}
\end{figure}

For the first numerical experiment we consider the right-hand side function (see Figure \ref{F1}: 
$$f_1(x,y)=\begin{cases}
x^2+y^2-1&\hbox{if }x^2+y^2<11,\\
\ds\frac{10}{1+\eps(x^2+y^2)^3}&\hbox{if }x^2+y^2>11.
\end{cases}$$

When $\ep$ and $\alpha$ are close to zero, the problem is then an approximation of the shape optimization problem
$$\min\left\{\int_\O u\,dx\ :\ -\Delta u=(x^2+y^2-1)\wedge 10\hbox{ in }\O,\ u\in H,\ u=0\hbox{ q.e. in }\R^d\setminus\O,\ |\O|\le m\right\}.$$
The solution of this last problem can be explicitly obtained and it is given by $\O=B(0,R)$ where $R$ is the biggest positive number satisfying $|\O|\leq m$ and such that the solution $u$ of the state equation is nonpositive, i.e.
$$R=\sqrt{m/\pi}\wedge \sqrt{2}.$$
In particular, we remark that the volume restriction is not saturated for $m>2\pi$. We see that the numerical solution \eqref{optpbnum} for the functions $g$ and $f_1$ given above is close to this one.

\begin{figure}[!]
\begin{minipage}[!t]{14cm}
\centering
\hspace{-.5cm}
\includegraphics[width=5.7cm]{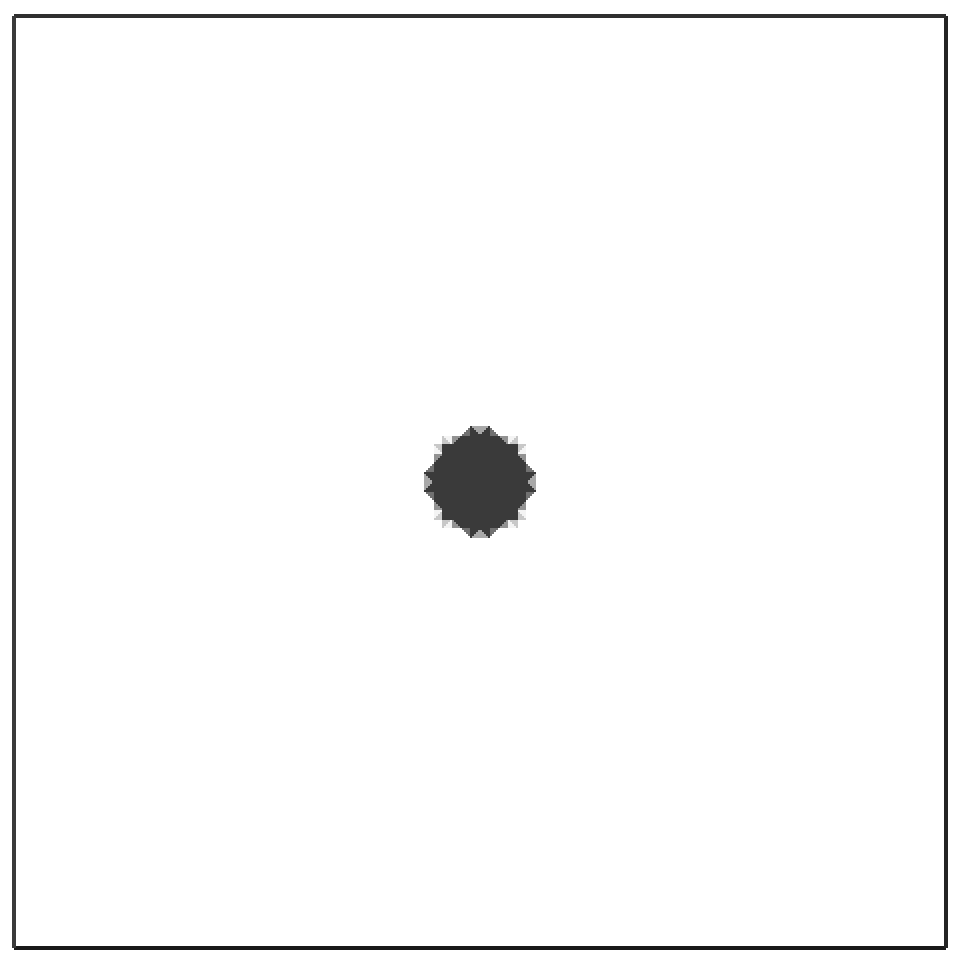}
\hspace{.7cm}
\includegraphics[width=5.72cm]{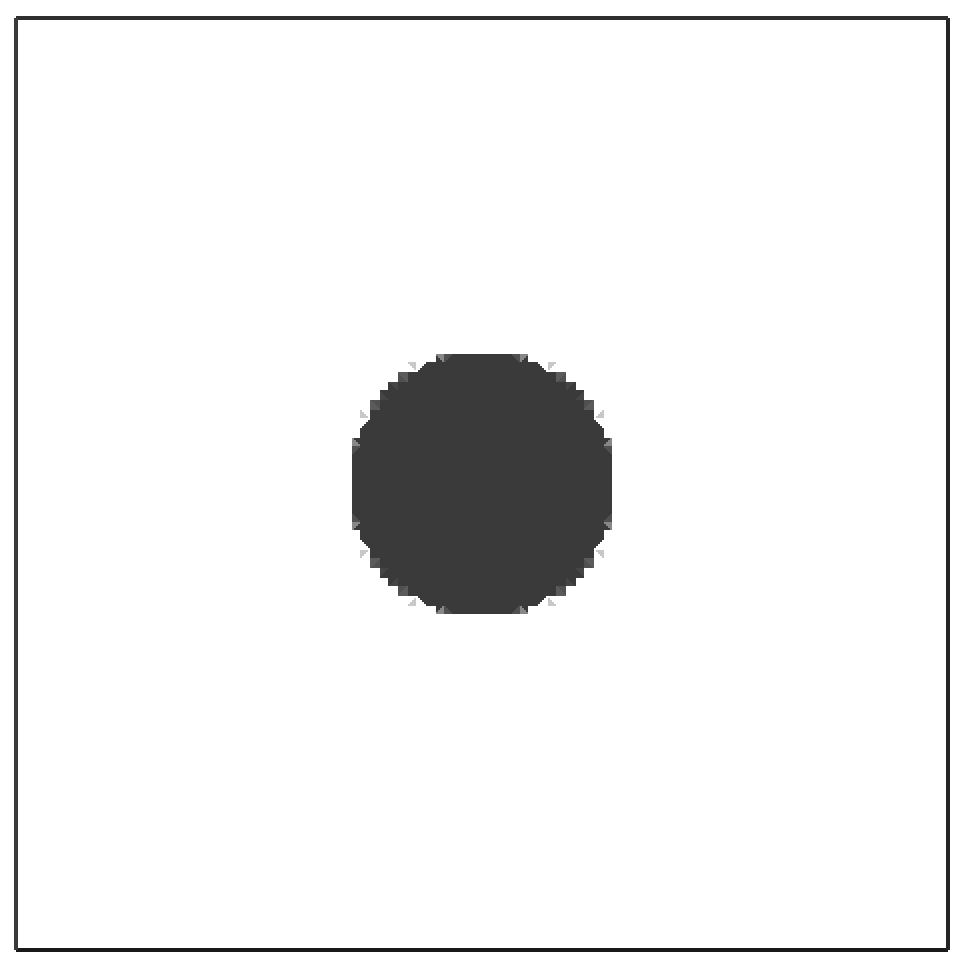}
\end{minipage}\vskip-0.1cm
\caption{The optimal potential $\mu_{opt}$ for volume constraint $m=2=m_{opt}$ (left) and $m=20> 6.367=m_{opt}$(right).} \label{F2}
\end{figure}
In Figure \ref{F2}, we show the optimal control provided by the numerical approximation for different constraints with respect to different amounts of available material. The black color represents the value $0$ for the optimal potential and the white color the value $+\infty$ (numerically $\mu_{max}$). We emphasize that in the obtaining this result we have not used the fact that the solution is known to be radial. As expected, in the left picture, corresponding to $m=2$, we can observe that the whole of the available material is placed in the part of the domain where the function $f$ is negative. Moreover, the volume constraint is saturated. In the right picture the amount of available material $m=20$ is bigger than the measure of the set where $f$ is negative. In this case, we can observe that optimal layout fulfills the set where $f$ is negative and also occupies a certain area where $f$ is positive. Now, the volume constraint is not saturated. In fact the amount of the material corresponding to optimal layout we find is $6.367$, which is close to the optimal value $2\pi\sim 6.283$ for the limit problem. We observe numerically that the optimum is independent of the choice of $M$, which suggests that effectively $\mu^s=0$ and $\mu^a$ are compactly supported.

\begin{figure}[!]
\begin{minipage}[!]{17cm}
\centering
\hspace{-1.9cm}
\includegraphics[width=7cm]{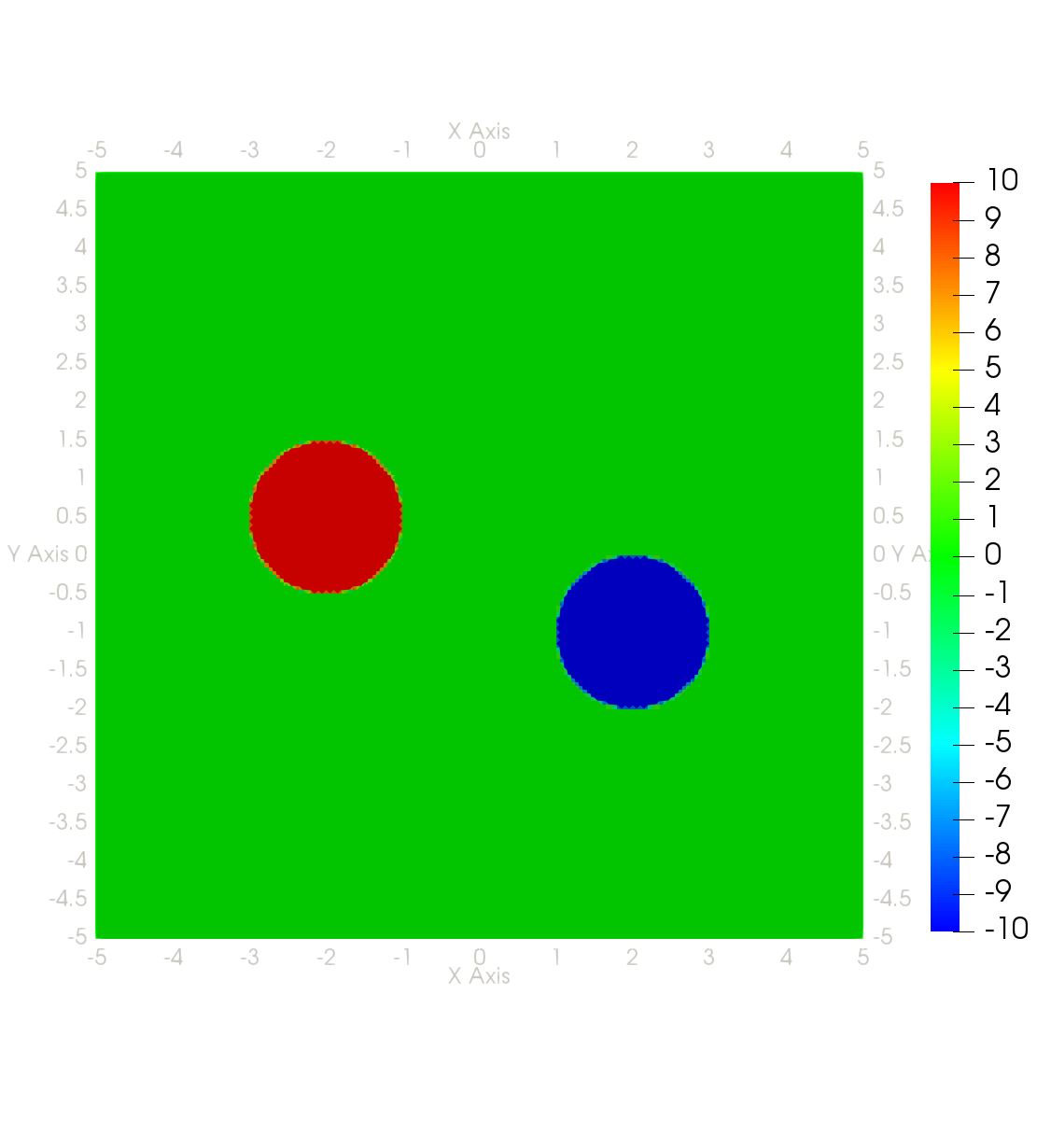}
\includegraphics[width=7cm]{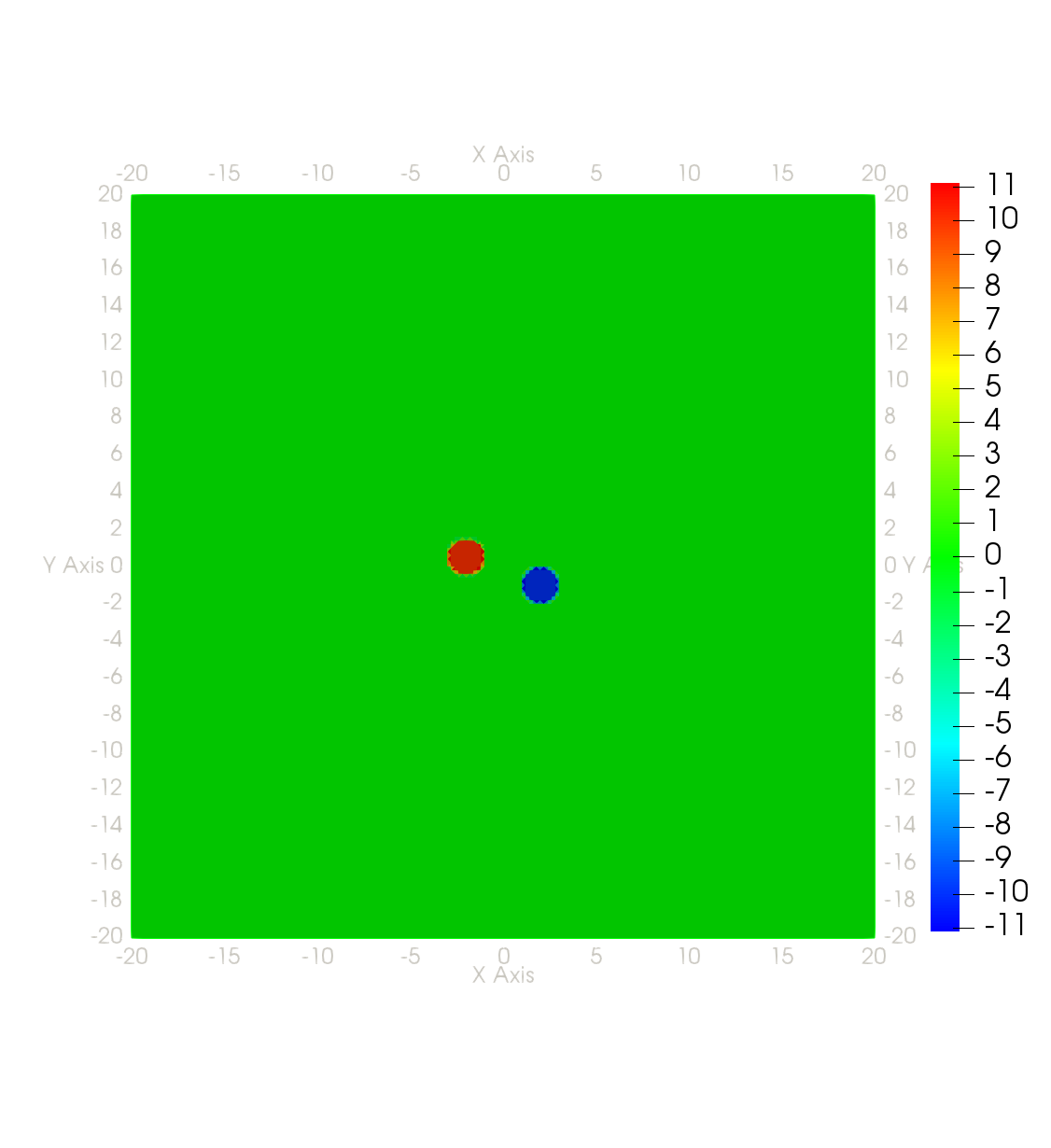}
\vskip-0.9cm
\end{minipage}\vskip-0.9cm
\caption{The rigth-hand side function $f_2$ in $D=(-5,5)\times(-5,5)$ (left) and in $D=(-20,20)\times(-20,20)$ (right).} \label{F3}
\end{figure}

\begin{figure}[!]
\begin{minipage}[!t]{14cm}
\centering
\hspace{-.5cm}
\includegraphics[width=6.9cm]{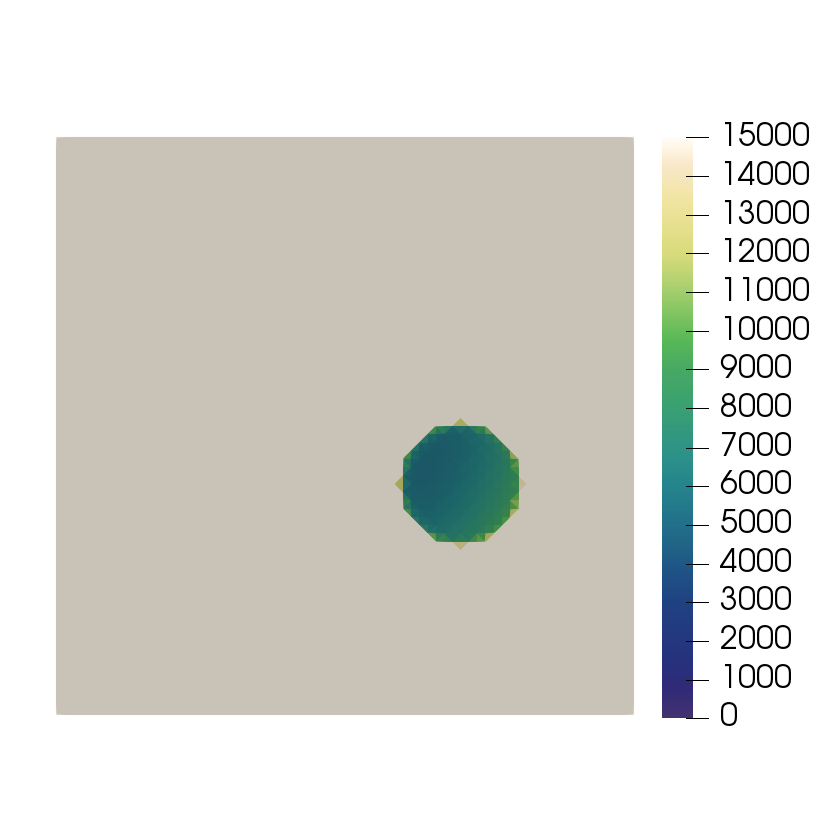}
\includegraphics[width=6.9cm]{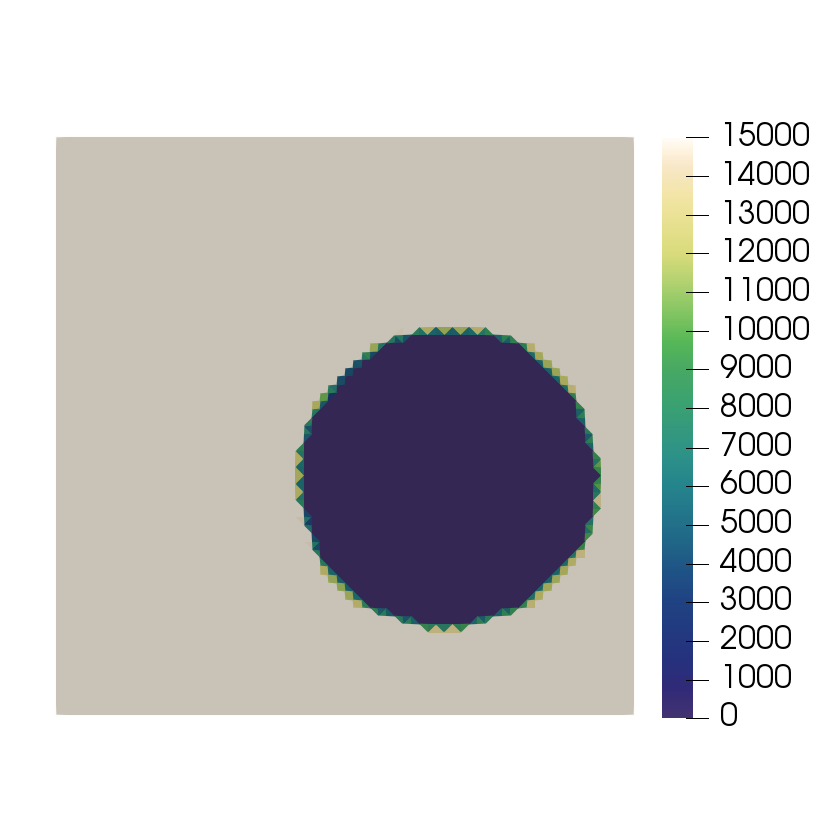}
\vskip-0.9cm
\end{minipage}\vskip-0.9cm
\caption{The optimal potential $\mu_{opt}$ in case 1 and $f_2$, $D=(-5,5)\times(-5,5)$, volume constraint $m=0.2$ (left) and $m=10$ (right).} \label{F6}
\end{figure}

\medskip
For the next example, we consider the right-hand side function:
\be\label{eq:f2}
f_2(x,y)=\begin{cases}
-10&\hbox{if }(x-2)^2+(y+1)^2<1,\\
10&\hbox{if }(x+2)^2+(y-0.5)^2<1\\
0&\hbox{otherwise}.
\end{cases}\ee
plotted in Figure \ref{F3}. Now, the solution cannot be radial and the optimal solution (even when $\ep$ and $\delta$ tend to zero) is not known. In Figure \ref{F6}, we have represented the optimal potential $\mu_{opt}$ corresponding to the domain $D=(-5,5)\times(-5,5)$ with different amounts of material. In the picture on the left, we consider the case where there is little material available, that is $m=0.2$. Then the numerical solution places the material where $f_2$ is negative and the volume constraint is saturated. In the right, we have consider a greater amount of material, that is $m=10$. In this case the optimal solution is a circular shape concentric with the disc where $f$ is negative and containing it. The volume constraint is saturated again.

\begin{figure}[!]
\begin{minipage}[!t]{14cm}
\centering
\hspace{-.5cm}
\includegraphics[width=6.9cm]{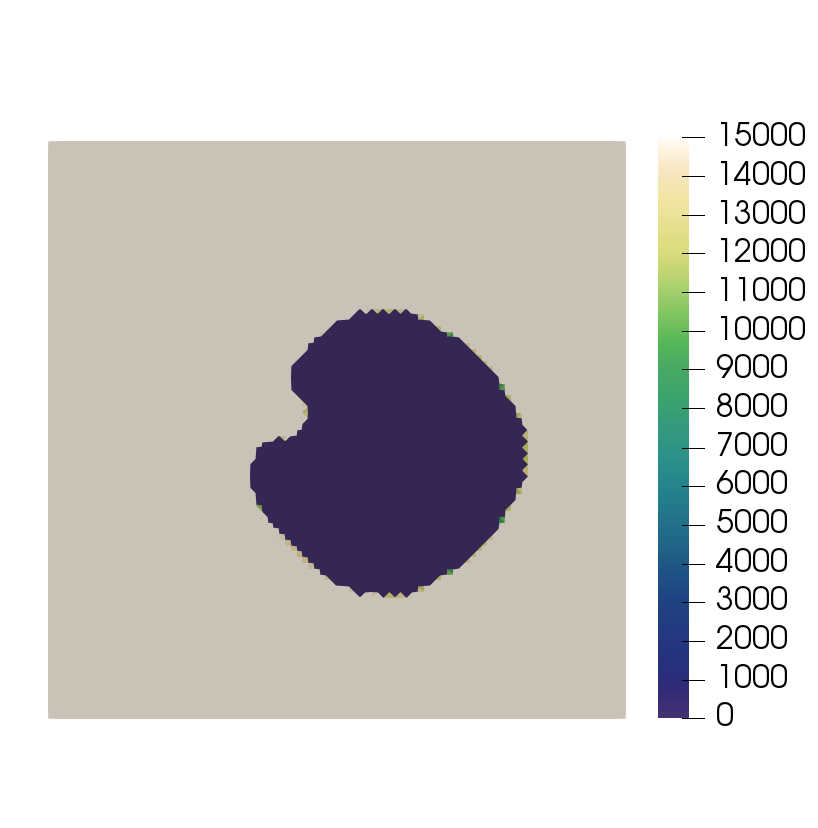}
\includegraphics[width=6.9cm]{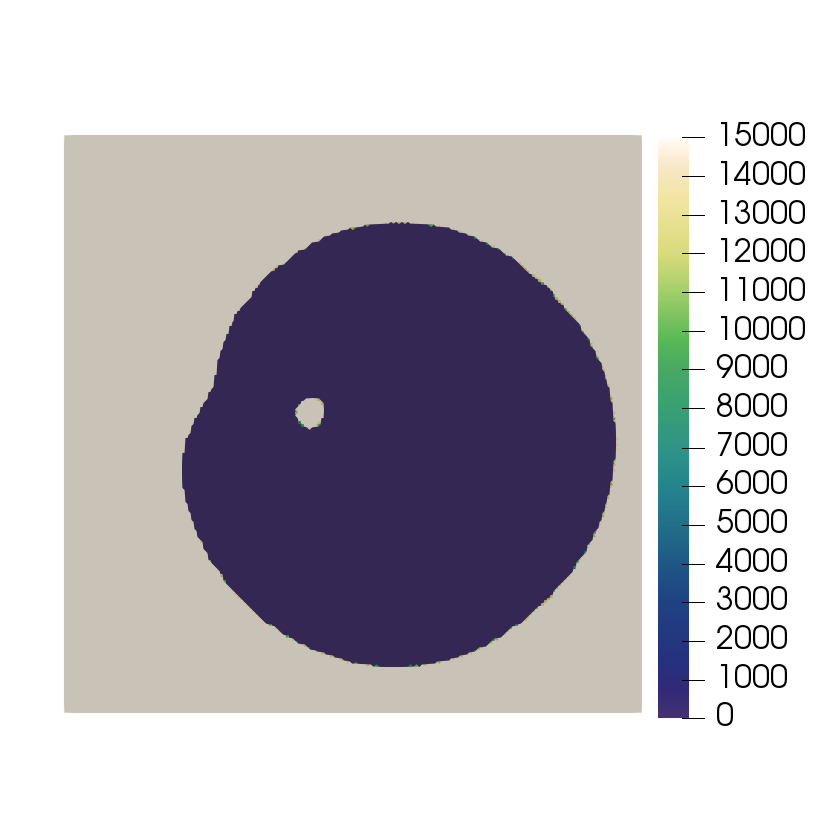}
\vskip-0.9cm
\end{minipage}\vskip-0.7cm
\caption{The optimal potential $\mu_{opt}$ in case 1 and $f_2$. $D=(-12.5,12.5)\times(-12.5,12.5)$ and volume constraint $m=110$ (left). $D=(-20,20)\times(-20,20)$ and $m=400$ (right).} \label{F7}
\end{figure}

In Figure \ref{F7}, we represent the evolution of the optimal solution where the amount of material increases. On the left we consider $m=110$ and $D=(-12.5,12.5)\times(-12.5,12.5)$. In this situation we observe that the optimal solution is again a characteristic function, but now it is not circular, because it tries to avoid the set where $f_2$ is positive. In the picture on the right we take $m=400$ and a bigger domain $D=(-20,20)\times(-20,20)$. Now we can see that the optimal shape is a quasi-circular region with a hole inside, corresponding with the area where $f_2$ is positive.

{\bf Case 2:} $\psi(s)=s^2$. For our simulations, we consider the right-hand side function $f_2$ defined in \eqref{eq:f2}. We observe that the volume restriction $\Psi(\mu)\le1$ is very different to the previous one. In the first case, to take $\mu=\infty$ does not require any expenses, while here it has an infinity cost and it is $\mu=0$ which is free.

\begin{figure}[!]
\begin{minipage}[!t]{14cm}
\centering
\includegraphics[width=6.9cm]{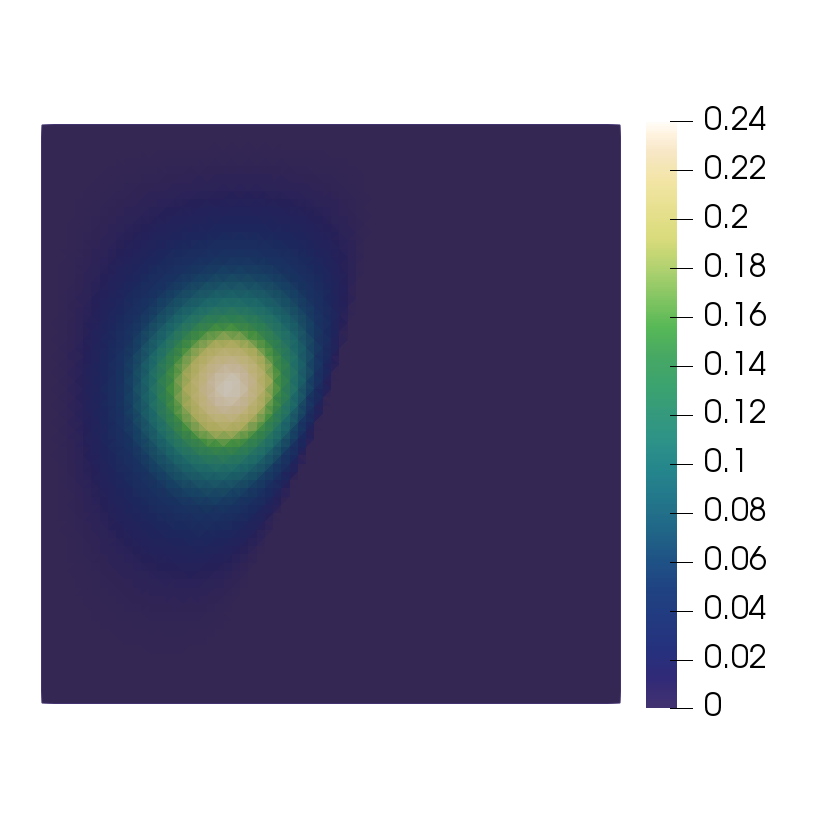}
\includegraphics[width=6.9cm]{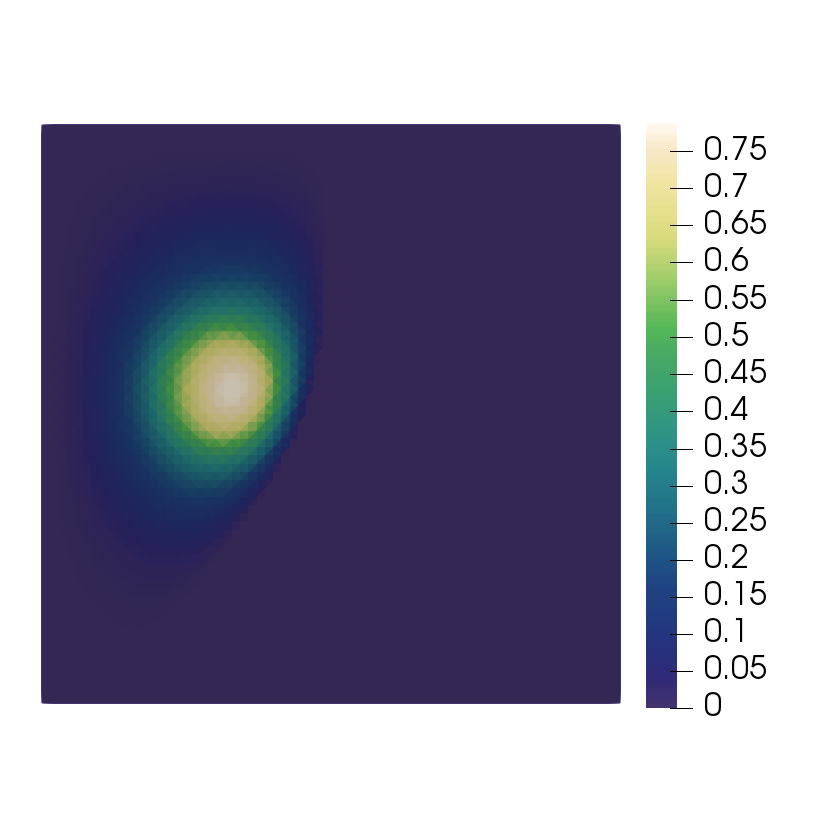}\\
\end{minipage}\vskip-0.1cm
\caption{The optimal potential $V_{opt}$ in case 2, $D=(-5,5)\times(-5,5)$, volume constraint $m=0.2$ (left) and $m=2$ (right).} \label{F4}
\end{figure}

\begin{figure}[!]
\begin{minipage}[!t]{14cm}
\centering
\includegraphics[width=6.9cm]{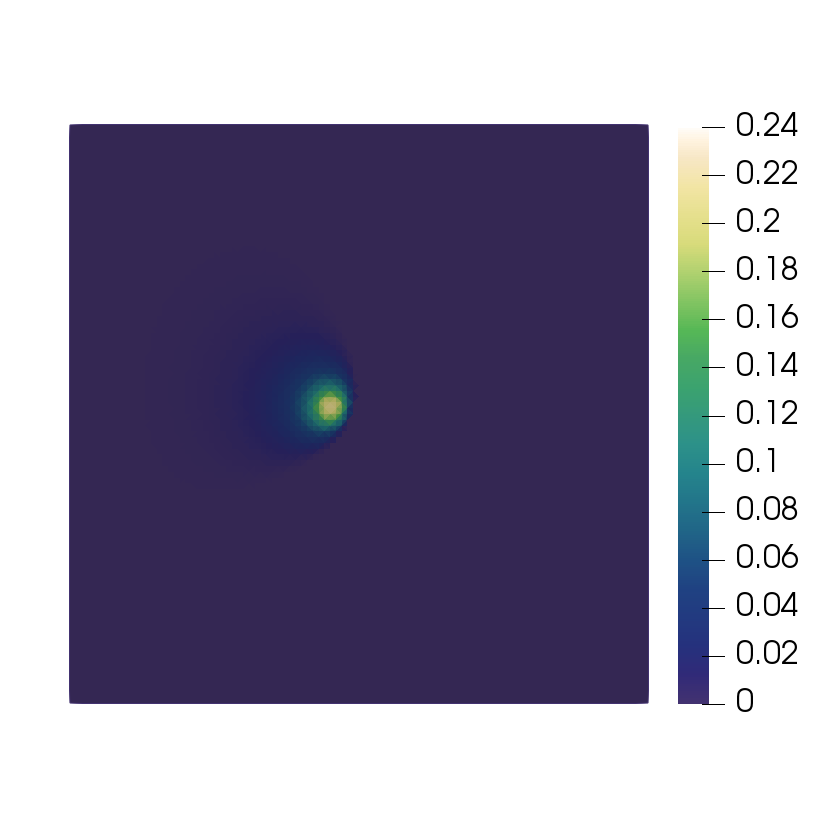}
\includegraphics[width=6.9cm]{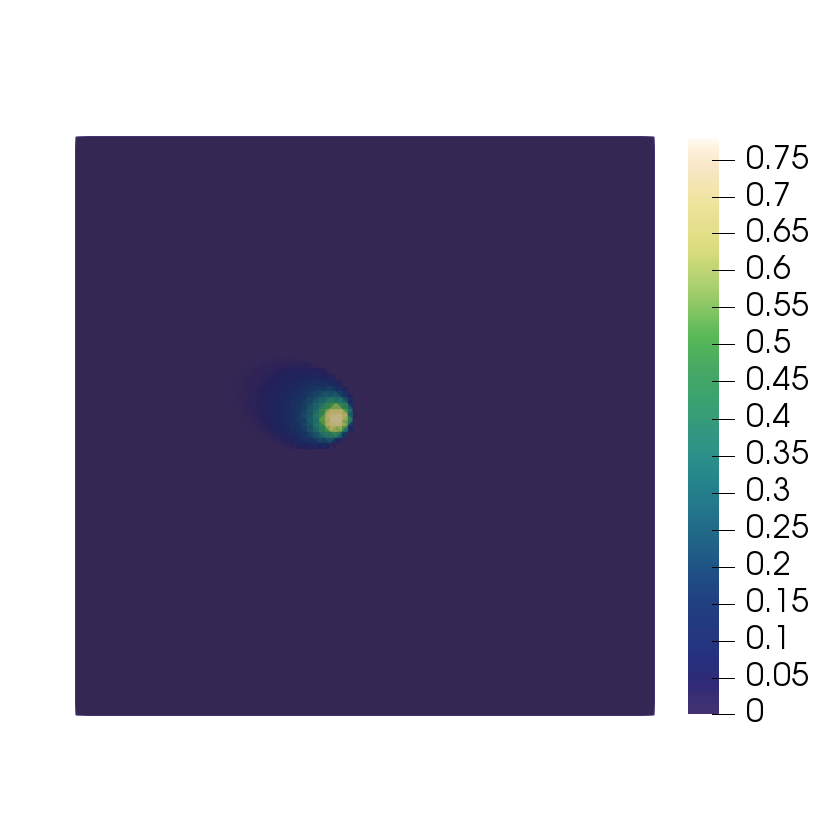}
\end{minipage}\vskip-0.1cm
\caption{The optimal potential $V_{opt}$ in case 2, $D=(-20,20)\times(-20,20)$, volume constraint $m=0.2$ (left) and $m=2$ (right).} \label{F5}
\end{figure}
Figures \ref{F4} and \ref{F5} show that the optimal potential $\mu_{opt}$ is finite in $D$. In Figure \ref{F4} we have solved the minimization problem in the domain $D=(-5,5)\times(-5,5)$ with two different values of available amount of material corresponding to $0.2$ on the left picture and $2$ on the right picture. In both cases the volume constraint is saturated, $\mu_{opt}$ is zero in most of $D$ and it is strict positive on the set where $f_2$ is positive. In Figure \ref{F5} we solve the problem in the domain $D=(-20,20)\times(-20,20)$. We observe that the corresponding solution in Figure \ref{F5} is just the extension by $0$ of the solution in Figure \ref{F4}, which suggests that $\mu^a$ is a bounded function with compact support.

\section*{Acknowledgements}
This work started during a visit of the first author at the Department of {\it Ecuaciones Diferenciales y An\'alisis Num\'erico} of the University of Sevilla and continued during a stay of the second and third authors at the Department of Mathematics of University of Pisa. The authors gratefully acknowledge both Institutions for the excellent working atmosphere provided. The work of the first author is part of the project 2017TEXA3H {\it``Gradient flows, Optimal Transport and Metric Measure Structures''} funded by the Italian Ministry of Research and University. The first author is member of the Gruppo Nazionale per l'Analisi Matematica, la Probabilit\`a e le loro Applicazioni (GNAMPA) of the Istituto Nazionale di Alta Matematica (INdAM). The second author and third authors have been partially supported by FEDER and ``Ministerio de Econom\'ia y Competitividad'' (Spain) under grant MTM2017-83583-P. 


\end{document}